\documentclass[11pt]{amsart}
 
\usepackage{times}
\usepackage[T1]{fontenc}
\usepackage{a4}
\usepackage{setspace}
\usepackage{amssymb}
\makeatletter
 \theoremstyle{plain}
 \newtheorem{thm}{Theorem}[section]
 \theoremstyle{plain}
  \newtheorem{cor}{Corollary}[section]
 \numberwithin{equation}{section} 
 \numberwithin{figure}{section} 
\theoremstyle{definition}
 \newtheorem{defn}[thm]{Definition}
 \theoremstyle{plain}
 \newtheorem{fact}{Fact}
 \theoremstyle{plain}
 \newtheorem{lem}[thm]{Lemma} 
 \theoremstyle{plain}
 \newtheorem{prop}[thm]{Proposition} 

\usepackage{psfrag}
 \usepackage{bbm}

\DeclareMathOperator{\Fix}{Fix}

\renewcommand{\Im}{\mathfrak{Im}}

\newcommand{\with}{\;|\;}

\newcommand{\R}{\mathbb {R}}

\newcommand{\C}{\mathbb {C}}
\newcommand{\1}{\boldsymbol{1}}

\newcommand{\N}{\mathbb {N}}

\newcommand{\e}{{\mathrm e}}
 \renewcommand{\i}{{\mathrm i}}

\newcommand{\w}{w}
\renewcommand{\v}{v}
\newcommand{\s}{\sigma}
\newcommand{\OA}{\mathcal{O}_A}
\newcommand{\FA}{\mathcal{F}_A}
\newcommand{\COA}{\mathcal{C}\left(\Sigma_A\right)}

\makeatother
\begin{document}
\title[Radon--Nikodym representations of Cuntz--Krieger algebras]{
Radon--Nikodym representations of Cuntz--Krieger algebras and
Lyapunov spectra for KMS states}
\author{Marc Kesseb\"ohmer}
\address{Fachbereich 3 - Mathematik und Informatik, Universit\"at Bremen, Bibliothekstrasse 1, 28359 Bremen, Germany}


\email{mhk@math.uni-bremen.de}

\author{Manuel Stadlbauer} \address{Institut f\"ur Mathematische Stochastik Maschm\"uhlenweg 8 - 10,  37 073 G\"ottingen, Germany} \email{stadelba@math.uni-goettingen.de}


 \author{Bernd O. Stratmann} \address{Mathematical Institute, University of St. Andrews, North Haugh,
 St Andrews,
 KY 16 9SS, Scotland} \email{bos@st-andrews.ac.uk}

\thanks{The second author was supported by the DFG project ``Ergodentheoretische Methoden in der hyperbo\-li\-schen Geometrie''.}

\date{January 14, 2006}
\begin{abstract}
We study  relations between $(H,\beta)$--KMS states on Cuntz--Krieger algebras
   and   the dual of the Perron--Frobenius 
operator $\mathcal{L}_{-\beta H}^{*}$.  Generalising  the well--studied purely hyperbolic 
situation, we obtain under mild conditions that for an expansive dynamical system there is a 
 one--one correspondence between $(H,\beta)$--KMS states and  eigenmeasures of 
$\mathcal{L}_{-\beta H}^{*}$ for the eigenvalue $1$. We then consider representations 
of Cuntz--Krieger algebras
which are induced by Markov fibred systems, and show that if the associated
incidence matrix is irreducible then these are $\ast$--isomorphic to the
given Cuntz--Krieger algebra. Finally, we apply 
these general results to study 
multifractal decompositions of  limit sets of essentially free Kleinian 
groups $G$ which may have parabolic elements. We show that 
for the Cuntz--Krieger algebra
arising from $G$ there exists an analytic family of 
KMS states   induced by
the Lyapunov spectrum of the analogue of the Bowen--Series map 
associated with $G$. Furthermore, we obtain a formula 
for the Hausdorff dimensions of the restrictions of these KMS states 
to the set of continuous functions on the limit set of $G$. If $G$ has 
no parabolic elements, then this formula 
can be interpreted as the singularity spectrum of the  
measure of maximal entropy associated with $G$.  
\end{abstract}
 \maketitle

\section{Introduction}
In mathematics and physics there is a long tradition of studying $(H,\beta)$--KMS 
(Kubo--Martin--Schwinger) states on various types of $C^{*}$--algebras
of observables. 
Originally, these notions stem from   quantum statistical mechanics, 
where $H$ refers to a given self--adjoint 
potential function (the {\em energy}) which fixes the system quantum 
mechanically, and $\beta$ admits the interpretation as 
the inverse of the {\em  temperature} of the system. The general
philosophy  here is that macroscopic thermodynamical properties
are reflected
within equilibrium states depending on $\beta$, whereas 
 microscopic quantum mechanical behaviour of the system 
is described by the $C^{*}$--algebra in combination with
some time evolution, that is 
a {\em gauge action} given by a one-parameter family 
$\left(\alpha_{H}^{t}\right)_{t \in \R}$ of 
$\ast$--automorphisms depending on $H$.

The first goal of this paper is to give a thorough review of the 
correspondence between fixed points of the dual of the Perron--Frobenius 
operator $\mathcal{L}_{-\beta H}^{*}$
and $(H,\beta)$--KMS states on a Cuntz--Krieger algebra $\OA$ 
associated with an  incidence matrix $A$. Under mild 
conditions on $\OA$, we collect 
various facts which mainly clarify the relation between the $(H,\beta)$--KMS states and the gauge action 
$\alpha_{H}^{t}$. The novelty here is that our approach includes the 
case in which the underlying dynamical system is expansive, and 
it therefore generalises the well--studied purely hyperbolic situation  
to the cases in which parabolic elements can occur.
More precisely, there is the well known result of Kerr/Pinzari (\cite{KerrPinzari:02}) and Kumjian/Renault (\cite{KumijanRenault:04}) which asserts that if $A$ is irreducible $H > 0$  then  there is a one--one 
correspondence between the eigenmeasures of 
$\mathcal{L}_{-\beta H}^{*}$ for the eigenvalue $1$
and $(H,\beta)$--KMS states.
We extend this result, using gauge--invariance of KMS-states (see Fact 8), by showing that if $H\geq 0$ then there exists a bijection from the set of eigenmeasures to the set of $(H,\beta)$--KMS states  which are trivial on $\{H=0\}$ (see Fact \ref{kms_gibbs}).

In Section 3 we then consider Markov fibred systems $(\Omega, m, 
\theta, \alpha)$ with an associated incidence matrix $A$.  We show that
if $A$ is irreducible and the Radon--Nikodym derivative of the 
measure $m$ is continuous on the support of $m$,
then  the associated Cuntz--Krieger 
$C^{*}$--algebra $\OA$  admits a representation 
 in terms of the Radon--Nikodym derivative of $m$  (see Theorem \ref{RN representation}).
 We will refer to this 
 representation as the Radon--Nikodym representation of $\OA$ induced 
 by the Markov fibred system. \\
 Finally, in Section 4 we give a finer 
 thermodynamical 
 ana\-ly\-sis of KMS states associated with the Cuntz--Krieger algebra $\OA$
 arising from an essentially free Kleinian group $G$.
More precisely, the final goal of this paper is to apply  the above general results 
to complete oriented 
$(n+1)$--dimensional hyperbolic mani\-folds whose
finitely generated non-elementary fundamental group 
$G$ is essentially free, that is $G$ is allowed to have parabolic 
elements and $G$ has no relations other than those ari\-sing
 from these parabolic elements.   
 The limit set $L(G)$ of $G$ is the smallest closed 
 $G$--invariant subset of the boundary ${\mathbb{S}}^{n}$ of 
 $(n+1)$--dimensional hyperbolic 
 space $ {\mathbb{D}}^{n+1}$. This set represents an intricate fractal set hosting the 
 complicated dyna\-mical action of $G$ on ${\mathbb{S}}^{n}$. 
 The combinatoric of this action 
 gives rise to some incidence matrix $A$, which then allows to 
 represent the 
 topological dynamics of the action  by a subshift of finite type $(\Sigma_{A}, \theta)$. In 
 particular, there exists  a topologically mixing Markov map 
 $T:L_{r}(G) \to L_{r}(G)$ on the complement $L_{r}(G)$ in $L(G)$
of the parabolic fixed points of $G$,  such that $T$  
 commutes with $\theta$ under the canonical coding map $\pi : 
 (\Sigma_{A},\theta) \to (L(G), T)$ (see Section 4 for the details). 
Functional analytic properties  of the action of 
$G$ on $L(G)$ can be studied by means
of  the  Cuntz--Krieger algebra $\OA$ associated with $G$,
and the link between $L(G)$ and $\mathcal{O}_{A}$ is given by the
commutative algebra
$\mathcal{C}\left(L(G)\right)$ of  complex-valued continuous functions on
$L(G)$. Note that $\OA$ is the
universal $C^{*}$--algebra obtained by completing the algebraic 
cross--product 
$\mathcal{C}(L(G)) \rtimes G $ with respect to the supremum of the norms 
of all $\ast$--representations on the underlying separable Hilbert 
space. It is well--known that in our situation here $\OA$ is nuclear, 
simple and purely infinite (\cite{A1}, \cite{A2}), and hence admits a 
classification by KK--theory (\cite{Ro}).  

 By considering a 
 special class $\{ I_{s}\}$ of potential functions, we show that there 
 exists a  family $\{\sigma_{s(\alpha)}\with \alpha \in 
 (\alpha_{-},\alpha_{+}) \}$ of $(I_{s(\alpha)},1)$--KMS 
 states on $\OA$, for some real analytic function $s: (\alpha_{-},\alpha_{+}) \to {\mathbb{R}}$,
 such that (see Theorem \ref{main}) 
 \[ \sigma_{s(\alpha)}(\log|T'|)=  \alpha,\]
 and such that for the Hausdorff dimension of the $\theta$-invariant weak
 $(-I_{s(\alpha)})$--Gibbs 
measure $\s_{s(\alpha)}\arrowvert_{\mathcal{C}(L(G))}$ we have 
\[\dim_H \left(\s_{s(\alpha)}\arrowvert_{\mathcal{C}(L(G))} \right)  
= \alpha^{-1} \cdot \s_{s(\alpha)}(I_{s(\alpha)} ). \]
Roughly speaking, these results are obtained by computing the Lyapunov 
spectrum of the map $T$. The proof in the convex cocompact case, that is if $G$ 
has no parabolic elements,  shows that these results can also be 
interpreted in terms of multifractal analysis as the multifractal 
spectrum or singularity spectrum of the  measure of maximal 
entropy arising 
from $T$. The proof for the parabolic case 
employs the method of inducing and here we follow closely our 
investigations of the dimension spectra 
arising from homological growth rates given in \cite{KesseboehmerStratmann:04a}. 
We remark that these results are based on a slight change of the 
usual approach which exclusively  considers multiplicative pertubations 
of the geometric potential 
function $J:=\log |T'|$. 
For this usual approach it is well--known 
that  in the convex cocompact case the exponent of convergence $\delta(G)$ 
is the unique $\beta$ for which there exists a $(\beta J,1)$--KMS 
state. However, if there are parabolic elements then additional to 
the Patterson-KMS state at $\delta(G)$  there exist  further 
$(\beta J,1)$-KMS states for each $\beta> \delta(G)$ which all correspond  to 
purely atomic measures concentrated on the fixed points of the 
parabolic transformations (see \cite{Sul1}, \cite{Sul2}, and for related 
investigations in terms of $C^{*}$--algebras for expanding dynamical systems  see  
 e.g. \cite{OP}, \cite{Evans}, \cite {R}, \cite{KumijanRenault:04}, 
 \cite{Exel}, \cite{MWY}, \cite{KerrPinzari:02}, \cite{Lott}).
In contrast to this usual approach,   in this paper
   we also allow additive pertubations of the potential functions.  
   More precisely,  our analysis is based on potentials of the
 form $sJ+P(-sJ)$, where $P(-sJ)$ refers to the topological pressure
 of the function $-sJ$. Note that if $G$ has parabolic elements then, in 
 order to get the thermodynamical formalism to 
 work,  we require $T$--invariant versions of the measures 
 derived from  the restrictions of the KMS states to 
 $\mathcal{C}(L(G))$.  This is subject of the appendix in which we 
 discuss a well--known formula of 
 Kac in the context of Markov fibred systems.
Recall that the Kac formula represents a convenient method which allows to derive an invariant 
measure  $\nu$
for the whole system  from a given measure $\widetilde{\nu}$ which is invariant under 
some induced 
transformation. The appendix gives a refinement of this well--known method 
by showing that there exist explicit 
formulae  which allow to compute the Radon--Nikodym derivative of $\nu$ 
with respect to the original transformation of the system in terms of the
Radon--Nikodym derivative of $\widetilde{\nu}$ with respect to the  
induced
transformation. These formulae might also be of  interest 
in general infinite ergodic theory.

\section{Cuntz--Krieger algebras, subshifts and KMS states}
Let $I$ refer to a given finite alphabet of cardinality at least $2$,  
and let $A=\left(a_{ij}\right)_{i,j\in I}$
be a transition matrix with entries in $\{0,1\}$ such that in each row
and column there is at least one entry equal to one.  The Cuntz--Krieger algebra
$\mathcal{O}_{A}$
associated with $A$, as
introduced in \cite{CuntzKrieger:80}, is a $C^{*}$--algebra generated by partial 
isometries $\left\{ S_{i}\with i\in I\right\} $
of a separable complex Hilbert space. The algebra $\mathcal{O}_{A}$
is the universal $C^{*}$--algebra generated by
$\{ S_{i}\with i\in I\}$ and satisfying the relations
 \begin{equation}
\sum_{j\in I}S_{j}S_{j}^{*}=1\,\,\,\,\textrm{and } \, \, \, S_{i}^{*}S_{i}=\sum_{j\in I}a_{ij}
S_{j}S_{j}^{*}, \;\textrm{ for all } i\in I.\label{CKrelation}\end{equation}
Recall that $S$ is a partial isometry if and only if $S=SS^{*}S$. Furthermore,
since  $\mathcal{O}_{A}$ is universal with respect to the relations
in (\ref{CKrelation}), we have with $ \delta_{ij} $ referring
to the Kronecker symbol,
\begin{equation} S_{i}S_{i}^{*}S_{j}S_{j}^{*}=\delta_{ij}S_{i}S_{i}^{*}  \textrm{ and }
S_{i}^{*} S_{i}S_{j}=a_{ij} S_{j} , \;\textrm{ for all } i,j\in I. 
\label{eq:CKRealtions2} 
\end{equation}
  It is well--known that  the algebra $\mathcal{O}_{A}$ is
uniquely determined up to isomorphism by the relations
in (\ref{CKrelation}) (see \cite[Theorem
2.13]{CuntzKrieger:80}). Furthermore, if $A$ is irreducible then
 $\mathcal{O}_{A}$ is
simple, which means that each  closed two--sided ideal of
  $\mathcal{O}_{A}$ is trivial (see \cite[Theorem 2.14]{CuntzKrieger:80}).

From a dynamical point of view, the transition matrix $A$ gives rise to a
subshift of finite type $(\Sigma_{A},\theta)$, where
\begin{eqnarray*}
  \Sigma_{A} &:=& \left\{ (\w_{0}\w_{1}\w_{2}\ldots)\in I^{\N}\with a_{\w_{i}\w_{i+1}}
  =1\,\,\,\textrm{for all } i \in \N_{0} \right\}, \\
 \theta&:&\Sigma_{A} \rightarrow \Sigma_{A},(\w_{0}\w_{1}\ldots)
 \mapsto(\w_{1}\w_{2}\ldots).
\end{eqnarray*}
The space  $\Sigma_{A}$ is a compact metric
space with respect to the metric given by
$$  d\left((\w_{0}\w_{1}\ldots),(\v_{0}\v_{1}\ldots)\right):=2^{-\min\left\{
 i\with \w_{i}=\v_{i}\right\} } , \textrm{ for } (w_{0}w_{1}\ldots),
(v_{0}v_{1}\ldots )\in \Sigma_{A}. $$
Also, it is well-known that $\theta$ is uniformly expanding with respect
to this metric (see e.g. \cite{DenkerGrillenbergerSigmund:76}).
Let $\mathcal{W}:=\bigcup_{n=1}^{\infty}\mathcal{W}^{n}$ refer to
the set of finite, admissible words, for
 \[
\mathcal{W}^{n}:=\{(\w_{0}\w_{1}\ldots\w_{n-1})\in
 I^{n}\with a_{\w_{i}\w_{i+1}}=1\;\textrm{for }i=0,1,\ldots n-1\}.\]
By defining for each $n\in\N$ and  $\w=(\w_{0}\w_{1}\ldots\w_{n-1}) \in \mathcal{W}^n$,
  \[
[\w]:=\left\{ \left(\w'_{0}\w'_{1}\ldots\right)\in 
\Sigma_{A}\with \w_{i}=\w'_{i}\:\textrm{for }i=0,\ldots,n-1\right\} ,\]
 we have that  $\{[\w]\with \w\in\mathcal{W}^{n}\}$ is a partition of $\Sigma_{A}$
 consisting of
closed and open sets.
The link between $\Sigma_{A}$ and $\mathcal{O}_{A}$ is then given by the
commutative algebra
$\mathcal{C}\left(\Sigma_{A}\right)$ of continuous, complex-valued functions on
$\Sigma_{A}$. With  $\1_{[w]}$ referring to the indicator function, we
have that the
algebra generated by $\left\{ \1_{[\w]}\with\w\in\mathcal{W}\right\} $
is dense in $\mathcal{C}\left(\Sigma_{A}\right)$. Therefore, the assignment
 $\psi\left(\1_{[\w]}\right):=S_{\w}S_{\w}^{*}$ admits an extension
to an isomorphism $\psi$ from $\mathcal{C}\left(\Sigma_{A}\right)$
to the commutative subalgebra generated by
$\left\{ (S_{\w} S_{\w}^{*})\with \w\in \mathcal{W} \right\}$, where
  $S_{\w}:=S_{\w_{0}}  S_{\w_{1}} \cdots S_{\w_{n-1}}$ for
  $\w = (\w_0 \w_1 \ldots \w_{n-1})$. Also, combining this with the second
  relation in (\ref{CKrelation}), one immediately verifies that $S_{\w}^{*} S_{\w} =
  \psi\left(\1_{\theta^n[\w]}\right)$, for all $\w \in
  \mathcal{W}^n$ and $ n\in \N$.
For ease of notation, throughout we will not distinguish between $f\in
\mathcal{C}(\Sigma_A) $ and $\psi\left(f\right) \in \mathcal{O}_A$.

  The following proposition gives some of the most
  important basic rules for the calculus within
  $\mathcal{O}_{A}$. In here, $\FA$ refers to the algebra of all finite linear
  combinations of words in the generators $S_j, S_j^*$. Note that $\FA$ is norm dense
  in $\OA$. Moreover, $\tau_\w : \theta^n([\w])\to [\w]$ denotes the
  inverse branch of $\theta^n|_{[\w]}$, for $\w \in \mathcal{W}^n$ and $n\in \N$.
  \begin{prop}\label{rules}
   For the Cuntz--Krieger algebra  $\mathcal{O}_{A}$,
   the following holds.
   \begin{itemize}
   \item[{\rm (1)}] Each  $X\in \FA$ can be written as a linear combination
    of terms of the form $S_\v S_\w^*$, for $\v,\w \in \mathcal{W}$.
   \item[{\rm (2)}] $S_\v^* S_\w = \delta_{\v,\w} \cdot S_\v^* S_\v$, for all $\v,\w \in
   \mathcal{W}^n$ and $n\in \N$.
   \item[{\rm (3)}] $S_{\v_0\ldots \v_{n-1}}^* S_{\v_0\ldots \v_{n-1}}
   = S_{\v_{n-1}}^* S_{\v_{n-1}}$,  for all $(\v_0\ldots \v_{n-1})\in
   \mathcal{W}^{n}$ and $n\in \N$.
   \item[{\rm (4)}] $S_\v^*f S_\w = \delta_{\v,\w}  \cdot  f \circ \tau_\v$, 
   for all $\v,\w \in \mathcal{W}^n$, $ n\in \N$ and $f \in \mathcal{C}\left( 
   \Sigma_A \right)$. In particular, we hence have  $S_\v^*f S_\w \in 
   \mathcal{C}\left( \Sigma_A \right)$ and $S_\v^*f S_\v = S_\v^* 
   S_\v \cdot f \circ \tau_\v$.
   \item[{\rm (5)}]  $S_\v f S_\v^* = \1_{[\v]} \cdot f \circ \theta^n$,  for all $\v \in
   \mathcal{W}^n$ and
   $n \in \N$. In particular, we hence have  $S_\v f S_\v^* \in 
   \mathcal{C}\left( \Sigma_A \right)$ and $S_\v f S_\v^* = S_\v  S_\v^* 
   \cdot f \circ \theta^n$.
   \item[{\rm (6)}] For all $\v \in \mathcal{W}^n$, $n \in  \N$ and $f \in \COA$, the following holds.
   \begin{itemize}
  \item[{\rm (a)}] $ S_\v f = (f \circ  \theta^n) S_\v$ and $ f S_\v  = 
  S_\v  (f\circ \tau_\v \cdot \1_{\theta^n[\v]})$.
  \item[{\rm (b)}]
  $ f S^*_\v  = S^*_\v  (f\circ \theta^n)$ and $ S^*_\v f = 
   ( f \circ \tau_\v \cdot \1_{\theta^n[\v]}  ) S^*_\v$.
  \end{itemize}
   \end{itemize}
  \end{prop}

  \begin{proof}
  For (1), (2) and (3) we refer to \cite[Lemma 2.1, 2.2]{CuntzKrieger:80}. For the proof of (4), we consider 
  without loss
    of generality $\v$ $=$ $(\v_0\ldots \v_{n-1})$, $\w =$ $( \w_0\ldots \w_{n-1})$ $\in$ $\mathcal{W}^n$
  and $u =(u_0 \ldots u_{m-1}) \in \mathcal{W}^m$, for $n < m$. Clearly, by (2) we have 
  that
  if either $\v  \neq (u_{0}\ldots u_{n-1})$ or $\w \neq (u_{0}\ldots 
  u_{n-1})$, then
  $S_\v^* S_u S_u^* S_\w = 0$.
  Otherwise, that is if $\v=\w=(u_0\ldots u_{n-1})$, we obtain from (2)
  and (3),
  \begin{eqnarray*}
   S_\v^* S_u S_u^* S_\w &= & S_\v^* S_{u_{0}\ldots u_{n-1}} S_{u_{n}\ldots u_{m-1}} S_{u_{n}\ldots
   u_{m-1}}^* S_{u_{0}\ldots u_{n-1}}^* S_\w \\
&=&    S_{u_{n-1}}^*S_{u_{n-1}}  S_{u_{n}\ldots u_{m-1}} S_{u_{n}\ldots u_{m-1}}^*
S_{u_{n-1}}^*S_{u_{n-1}}\\
    &=&   \1_{\theta^n[u]} =  \1_{[u]}\circ \tau_\v .
  \end{eqnarray*}
  Since $\{ \1_\w \with \w \in \mathcal{W}\}$ is dense in $\mathcal{C}\left( \Sigma_A \right)$, the result
  in (4)
  follows. The assertion in (5) follows by similar means.

  For the proof of (6), we consider without loss of generality  $ \1_{[\w]} = S_\w S_\w^* \in \COA$,
  for $\w =(\w_0\ldots \w_{m-1})\in \mathcal{W}$.
  Since $S_\w S_\w^* S_j S_j^*  \1_{[\w]} \1_{[j]} = \delta_{\w_0 j}  \1_{[\w]}$,
  we have for  $\v=(\v_0\ldots \v_{n-1})$
  $\in \mathcal{W}^{n}$ with  $a_{\v_{n-1} \w_{0} } =1$, using (\ref{CKrelation}),  (3) and (5),
  \begin{eqnarray*}
  S_\v  S_\w S_\w^*
   = S_\v  S_\w S_\w^* \sum_{j \in I} a_{\v_{n-1},j} S_j S_j^* = S_\v  \1_{[\w]} S_v^*S_v = (\1_{[\w]}\circ
   \theta^n) S_v .
  \end{eqnarray*}
 This gives the first part in (a). For the second part note that
 \[ f S_\v = f  S_\v S^*_\v S_\v = f \cdot \1_{[\v]} S_\v =
  (f \circ \tau_\v \cdot \1_{\theta^n[\v]})\circ
 \theta^n S_\v. \]
Therefore, using the first part, the assertion follows. The statements in (b)
 are immediate consequences of (a).
  \end{proof}

Throughout, let  $H = H^* \in \COA$ be always a given self--adjoint potential function.
Then the \emph{Perron--Frobenius operator}
$\mathcal{L}_{H}:
\mathcal{C}\left(\Sigma_{A}\right)\to\mathcal{C}\left(\Sigma_{A}\right)$ 
associated with $H$ is defined in dynamical terms as follows (see \cite{Ruelle}).
 For $f \in \mathcal{C}\left(\Sigma_{A}\right) $ and $x\in \Sigma_A$, let
\[ \left(\mathcal{L}_{H}\left(f\right)\right)
\left(x\right):=\sum_{y\in \theta^{-1}\left(\left
\{ x\right\} \right)}e^{H\left(y\right)}f(y).\]
 Using Proposition \ref{rules} (4), we obtain that  $\mathcal{L}_{H}$ 
 can be expressed in algebraic terms in the following way.
\begin{equation}
   \mathcal{L}_{H}\left(f\right) =  \sum_{j \in I}S_{j}^{*}
e^{H}f S_{j}  .\label{dual}
\end{equation}

\begin{defn}\label{gauge-KMS} 
 For $ t,\beta \in \R$, we define the following.
   \begin{itemize}
  \item A \emph{gauge action}\/
is a $\ast$--automorphism $\alpha_{H}^{t}:\mathcal{O}_{A}\to\mathcal{O}_{A}$
given by the extension to $\OA$ of 
 \[ 
 S_{j} \mapsto \alpha_{H}^{t}S_{j} := e^{itH}S_j , \textrm{ for each } j \in I.
 \]
\item   A $\left(H,\beta \right)$\emph{--KMS
state} is a  state   $\s$ on $\mathcal{O}_{A}$ such that for each pair $X,Y$ in a norm--dense
subset of $\mathcal{O}_{A}$, there exists an analytic function $F_{X,Y} : 
\{ z \in \C \with 0 \leq  \Im (z) \leq \beta \} \to \C$
 such that
\[ F_{X,Y}(t) = \s\left(X \alpha_{H}^{t}(Y)\right)\textrm{ and }
 F_{X,Y}(t + i\beta) = \s\left( \alpha_{H}^{t}(Y) X\right).\]
  \end{itemize}
\end{defn}
Recall that a state $\s$ on a $C^*$--algebra $\mathcal{A}$ is by definition a linear
functional for which $\|\s \|=1 $ and $ \s(X^*X) \geq 0$, for all $X
\in \mathcal{A}$.
Also note that, since $ \alpha_{H}^{t}$ is a  $\ast$--automorphism, we have for all
$\v = (\v_0  \cdots \v_{m-1})$, $\w = (\w_0  \cdots \w_{n-1}) \in \mathcal{W}$, and
$t \in \R$,
 \begin{eqnarray}
\alpha_{H}^{t}\left(S_{\v}S_{\w}^{*}\right)&=& e^{\i {t} H}S_{\v_{0}}\cdots e^{\i
  {t}H}S_{\v_{m-1}}S_{\w_{n-1}}^{*}e^{-\i t H} \cdots S_{\w_{0}}^{*}e^{-\i
  t H} \label{eq:GaugeAktionOnGenerators}  \\
  &=& e^{\i t \sum_{k=0}^{m-1}H\circ \theta^k }S_\v S_\w^*
   e^{-\i t \sum_{k=0}^{n-1}H\circ \theta^k }.\nonumber
  \end{eqnarray}
In here, the final equality is a consequence of Proposition \ref{rules} (6).
In order to consider analytic continuations of the gauge action,
we require the following concept of analyticity of 
\cite[Definition 2.5.20]{BratteliRobinson:87}.

\begin{defn}\label{analytic}
An element $X \in \mathcal{O}_A$ is called \emph{$\alpha_H^t$--analytic} 
if there exists a positive number $\lambda$ and a map 
 $f_X : D_\lambda := \{ z \in \C \with |\Im(z)|<\lambda \} \to 
 \mathcal{O}_A$ such that the following holds.
\begin{enumerate}
\item $ f_X(t) = \alpha_H^t(X)$, for all $t \in \R$.
\item  The function $\eta \circ f_X : D_\lambda \to \C $ is
analytic, for each map $\eta$ in the topological dual $(\OA)'$ of $\OA$.
\end{enumerate}
\end{defn}

\begin{fact}[\textbf{$\alpha_{H}^{t}$-analyticity}] 
\label{Analytic} The gauge 
action admits a unique continuation to all of\/ $\C$ in the  
following way. For each $X=S_{\v}  S^*_{\w}$ with  $m,n \in \N$, $\v \in 
\mathcal{W}^{m}$ and $\w \in \mathcal{W}^{n}$,
there exists a unique continuation of $\alpha_{H}^t(X)$  such that 
\[  \alpha_{H}^{z} (X) = e^{\i z \sum_{k=0}^{m-1}H\circ \theta^k }  X 
   e^{-\i z \sum_{k=0}^{n-1}H \circ \theta^k}, \textrm{  for all  } z 
    \in \C. \]
    In particular, we hence have that each $Y \in \FA$ is 
    $\alpha_H^t$--analytic with respect to  $D_\lambda = \C$. 
\end{fact}

\begin{proof} 
Note that if $X \in \mathcal{O}_A$ is $\alpha_H^t$--analytic  and 
$\eta \in (\OA)'$, then analyticity of 
$\eta \circ f_X$ together with Definition \ref{analytic} (1)  immediately 
implies that the function $\eta \circ f_X$ does not depend on the special choice of $f_X$. 
Combining this with the fact that $\OA$ is a Hilbert space and hence is reflexive,  
the uniqueness of $f_X$ follows.  For the remaining assertions, 
first note that by \cite[Proposition 2.5.21]{BratteliRobinson:87} we have 
that the statement in Definition \ref{analytic} (2) is equivalent to the 
fact that  for each $z \in D_\lambda$ the following limit exists, where 
the limit is taken  with respect to the 
norm in $\OA$. 
\[ \lim_{w \to z} \frac{f_X(z)-f_X(w)}{z-w}. \]
Hence, in order to show that each $Y \in \FA$ is 
$\alpha_H^t$--analytic, it is sufficient to show that the above limit exists,  
for  $z \in \C$ and $X={S_\v  S^*_\w}$  with $\v \in \mathcal{W}^{m}$,
 $\w \in \mathcal{W}^{n}$ and 
\[ f_{S_\v  S^*_\w}(z) := e^{\i z \sum_{k=0}^{m-1}H\circ \theta^k }
 S_\v S_\w^* e^{-\i z \sum_{k=0}^{n-1}H\circ \theta^k}.
\]
Indeed, for $ g := \sum_{k=0}^{m-1}H\circ \theta^k - \sum_{k=0}^{n-1}H\circ 
\theta^k \circ \tau_\w \circ \theta^m$ we have by 
 Proposition \ref{rules} (6) that $ f_{S_\v  S^*_\w}(z) = e^{\i z g} 
 S_\v  S^*_\w$. Hence, it follows for $z, \varepsilon \in \C$,
\[
\lim_{\varepsilon \to 0} \frac{ f_{S_\v  S^*_\w}(z)-f_{S_\v  S^*_\w}
(z + \varepsilon)}{ \varepsilon}  
=   \left(\lim_{\varepsilon \to 0}  \frac{1- e^{\i \varepsilon g} }
{\varepsilon}\right) \cdot e^{\i z g} S_\v  S^*_\w 
=  -\i g   \cdot f_{S_\v  S^*_\w}(z).
\]
This shows that $S_\v  S^*_\w$ is $\alpha_H^t$--analytic with respect to  
$D_\lambda = \C$, which then clearly also holds for each $Y \in \FA$. 
Since, as we have seen above, $f_X(z)$ is uniquely determined for all 
$z \in \C$, we can now define  
$\alpha^z_H(X):=f_X(z)$, which then finishes the proof.   
\end{proof}

Note that for $z \in \C \setminus \R$ the continuation  $\alpha^z_H$
of the gauge action  $\alpha^t_{H}$ is no longer  a $*$--automorphism.

\begin{fact}[\textbf{KMS condition}]\label{KMS condition}
 For an $(H,\beta)$--KMS state $\s$ on $\mathcal{O}_{A}$ with $\beta \in \R
 \setminus \{0\}$, we have
 
\[  \s(XY) = \s(Y\alpha_H^{\i\beta}(X)) , \textrm{ for 
 all }  X,Y \in \FA.\]
\end{fact}

\begin{proof} 
  A similar argument as in the proof of Fact \ref{Analytic}  shows that
   the assignment $z \mapsto \s(Y \alpha^z_H(X))$ gives 
  rise to an analytic function on $\C$, for each $(H,\beta)$--KMS state 
  $\s$ and each $X,Y \in \FA$.  This implies 
  $\s(Y\alpha^z_H(X)) = F_{Y,X}(z)$, for all $ z \in \C$ 
  for which $ 0 \leq  \Im (z) \leq \beta $. Using Definition 
  \ref{gauge-KMS} (2), it therefore 
  follows that
\[\s(XY) = F_{Y,X}(\i\beta)  = \s (Y f_{X}(i\beta)) = 
\s(Y\alpha_H^{\i\beta}(X)) .\]
\end{proof}

We remark that the proof of the previous fact in particular  shows that
 \[  \s(\alpha_H^t(X)Y) = \s(Y\alpha_H^{t+ \i\beta}(X)) ,
  \textrm{ for  all } t\in \R   \textrm{ and }  X,Y \in \FA.\]

For a KMS state we also immediately obtain the following fact (see e.g. 
\cite[Proposition 5.3.3]{BratteliRobinson:97}).
\begin{fact}[\textbf{State invariance}]\label{betainvariance}
 For an $(H,\beta)$--KMS state $\s$ on $\mathcal{O}_{A}$ with $\beta \in \R
 \setminus \{0\}$, we have
 \[ \s\left(\alpha_{H}^{z}X\right)=\s\left(X\right), \textrm{ for all } 
 z\in \C   \textrm{ and } X \in \OA.\]
\end{fact}

Let us collect further important observations concerning the thermodynamical
formalism for KMS states. The following fact is adopted from
\cite[Section 7]{KerrPinzari:02}, and  we include a proof for convenience.

\begin{fact}[\textbf{Centraliser}] 
\label{fact-invariance} For an $(H,\beta)$--KMS state $\s$ on $\mathcal{O}_{A}$, we have
\[ \s(Xf) = \s(fX), \textrm{ for all }f \in \COA, X \in \OA. \]
\end{fact}

\begin{proof} Using (\ref{eq:GaugeAktionOnGenerators}) and Proposition
\ref{rules} (6), we obtain for all $\w \in\mathcal{W}$ and $t \in \R$,
\[
\alpha_{H}^{z}\left(S_{\w}S_{\w}^{*}\right)=  S_{\w}S_{\w}^{*}.
\]
Since $\left\{ S_{\w}S_{\w}^{*}: \w \in\mathcal{W}\right\} $ is dense in 
$\mathcal{C}\left(\Sigma_{A}\right)$, it
follows that $\alpha_H^z(f)=f$, for all $f \in \mathcal{C}\left(\Sigma_{A}\right)$. 
Using Fact \ref{KMS condition}, we hence have 
\[   \s(Xf) = \s(X\alpha_H^{\i \beta}(f)) = \s(fX) .\]
\end{proof}

\begin{fact}[\textbf{Faithfulness}] \label{faithful}
 Let $A$ be irreducible, and let $\s$ be an $(H,\beta)$--KMS state. We then have for all  
  $X \in \mathcal{O}_{A}$,
  \[ \s(X^*X)=0  \, \hbox{  if and only if  }  \, X=0.\]
\end{fact}
\begin{proof} First, note that the
set \[
\mathcal{I}:=\left\{ Y\in\mathcal{\mathcal{O}}_{A} \with \s\left(Y^{*}Y\right)=0\right\} \]
is a left ideal of $\mathcal{\mathcal{O}}_{A}$.  This follows since, using the estimate
$|\s(X^*YX)|\leq \s(X^*X)\|Y\|$ (see \cite[Proposition 2.3.11]{BratteliRobinson:87}),  
we have for $Y\in\mathcal{I}$ and $X\in\mathcal{\mathcal{O}}_{A}$, 
 \[
0\leq\s\left(\left(X Y\right)^{*}X Y\right)\leq\left\Vert X\right\Vert 
^{2}\s\left(Y^{*}Y\right)=0.\]
Also, since $\s$ is continuous, $\mathcal{I}$ is closed.
Using the KMS condition and the Cauchy--Schwarz inequality for $\s$ 
(see \cite[Lemma 2.3.10]{BratteliRobinson:87}
and proof of Fact \ref{gaugeinvariance}), we obtain
\begin{eqnarray*}
\left|\s\left(X^{*}Y^{*}YX\right)\right|^{2} & = & 
\left|\s\left(\alpha_{H}^{-\i \beta}\left(X\right)X^{*}Y^{*}Y\right)\right|^{2}\\
 & \leq & \s\left(\left(\alpha_{H}^{-\i \beta}\left(X\right)X^{*}\right)^{*}
 \alpha_{H}^{-\i\beta}
 \left(X\right)
 X^{*}\right)\s\left(Y^{*}Y Y^{*}Y\right)\\
 & \leq&
 \s\left(\left(\alpha_{H}^{-\i\beta}\left(X\right)X^{*}\right)^{*}\alpha_{H}^{-\i \beta}
 \left(X\right)X^{*}\right)\left\Vert Y\right\Vert ^{2}\s\left(Y^{*}Y\right)=0.\end{eqnarray*}
Since $\mathcal{\mathcal{O}}_{A}$ is simple  if $A$ is irreducible 
(see \cite[Theorem 2.14]{CuntzKrieger:80}) and
since $1\not\in\mathcal{I}$, it follows that $\mathcal{I}=0$.
\end{proof}

\begin{fact}[\textbf{$\beta$ - Conformality}] \label{conformality}
For an $(H,\beta)$--KMS state $\s$ on $\mathcal{O}_{A}$ and for all
 $\v \in \mathcal{W}^m$, $\w \in \mathcal{W}^n$ with $m,n \in \N$ such that $\v\w \in \mathcal{W}^{n+m}$ we have 
\begin{enumerate}
\item $ \displaystyle \s\left(S_{\v\w} S_{\v\w}^* \right) = \s\left(e^{-\beta \sum_{k=0}^{m-1}
 H\circ \theta^k \circ \tau_\v} S_{\w} S_{\w}^*\right) $,\\[-2mm]
\item $ \displaystyle  \s\left(S_{\w} S_{\w}^* \right) = \s\left(e^{\beta \sum_{k=0}^{m-1} 
H\circ \theta^k } S_{\v\w} S_{\v\w}^*\right).$
\end{enumerate} 
\end{fact}
\begin{proof} By the KMS property, Fact \ref{Analytic} and Proposition \ref{rules} (6), it follows that
\begin{eqnarray*}
\s\left(S_{\v\w} S_{\v\w}^* \right)  &=& \s\left(S_{\w} S_{\w}^* 
 S_{\v}^* \alpha_{H}^{\i \beta}\left(S_\v\right) \right) 
= \s\left(S_{\w} S_{\w}^*  S_{\v}^* e^{-\beta\sum_{k=0}^{m-1} H\circ \theta^k} S_\v \right) \\
&=& \s\left(S_{\w} S_{\w}^*  S_{\v}^* S_\v e^{-\beta\sum_{k=0}^{m-1} 
H\circ \theta^k \circ \tau_\v} \right) \\
&=& \s\left(S_{\w} S_{\w}^*  e^{-\beta\sum_{k=0}^{m-1} H\circ \theta^k \circ \tau_\v} \right).
\end{eqnarray*}
The second assertion follows by precisely the same means. \end{proof}

We remark that Fact \ref{conformality} has the following immediate 
implication for  the measure  associated with $\s$ by the Riesz representation 
theorem. Namely, with 
$m_\s:=\s|_{\mathcal{C}(\Sigma_A)}$ referring to this measure, the second 
statement in Fact \ref{conformality} immediately gives   
\[ \log \frac{d m_\s \circ \theta^n}{d m_{\s}} = \beta \sum_{k=0}^{n-1} 
H\circ \theta^k, \hbox{ for each } n \in \N.\]

For the following we recall the definitions of a Gibbs measure and  of a weak Gibbs measure,
which we have adapted to our situation here (see \cite{Kesseboehmer:01}, \cite{Yuri}). 
\begin{defn}\label{definition:weak Gibbs} Let $\mu$ a Borel probability measure on 
$\Sigma_{A}$ for which there exists $f \in \COA$ and   
a  sequence $\left(b_{n}\right)_{n \in \N} $ of positive numbers such that, for all $n \in \N$, $\w \in
  \mathcal{W}^n$ and $x \in [w]$,
    \[ e^{-b_{n}} \leq   \frac{\mu([w]) }{e^{\sum_{k=0}^{n-1} f \circ
  \theta^k (x)} } \leq e^{b_{n}} .\]
\begin{enumerate}
  \item The measure   $\mu$ is called  \emph{$f$--Gibbs measure} if the sequence\/ $(b_n)$  is  constant.
  \item The measure   $\mu$ is called \emph{weak $f$--Gibbs measure} if\/  $\lim_{n\to \infty}
  b_{n}/n=0$.
\end{enumerate}  
 \end{defn}
Recall that if $f$ is a
strictly negative  H\"older continuous potential 
function then the Perron--Frobenius--Ruelle theorem (see e.g. 
\cite{Bowen}) implies that 
there exists a unique $f$--Gibbs measure. Also, note that  the concept  
of a weak $f$--Gibbs measure 
is slightly more 
general than the concept of an  $f$--Gibbs measure. Namely, if for 
instance $f \in \mathcal{C}(\Sigma_{A})$ is not H\"older continuous  and $f\leq 0$ such that 
 $f(x) =0 $ for
at most finitely many $x \in \Sigma_A$, then it is possible that there exists a weak 
$f$--Gibbs measure which is 
not an $f$--Gibbs measure. In particular, this weak 
$f$--Gibbs measure  is not necessarily unique (see \cite{Kesseboehmer:01}, \cite{Yuri}). 
Next recall that the topological pressure of $f \in \COA$ is given by
\[ P(f) :=  \lim_{n \to \infty} \frac{1}{n}     \log  \sum_{\w \in \mathcal{W}^n} \exp
\left( \sup_{x \in [w]} \sum_{k=0}^{n-1} f\circ \theta^k(x)      \right) .\]
One immediately verifies that for a potential $f \in \COA$ for which there exists  a weak 
$f$--Gibbs measure, we necessarily have that  $P(f) =0$.
   For the following fact note that in  \cite{Kesseboehmer:01} it was
shown that if  $\mu$ is contained in the set $\Fix
\left(\mathcal{L}_{f}^{*}\right) $ of probability
measures $\nu$ which are eigenmeasures of the
dual $\mathcal{L}_{f}^{*}$ of the Perron--Frobenius
operator for the eigenvalue $1$, that is for
which $ \mathcal{L}_{f}^{*}\nu=\nu$,
then $\mu$ is a weak $f$--Gibbs measure. In this situation we then have
\[ \frac{d\mu \circ \theta}{d\mu} = e^{-f}.\]
Note that a weak $f$--Gibbs measure has no atoms if $f$
is strictly negative.
\begin{fact}[\textbf{$\mathcal{L}^*$--invariance}]\label{L*invariance}
 For an  $\left(H,\beta\right)$--KMS state $\s$ on
  $\mathcal{O}_{A}$ with $\beta \in \R \setminus \{0\}$, we have
\[ \sum_{j\in I}\s\left(S_{j}^{*}\e^{-\beta H}XS_{j}\right)=\s\left(X\right), 
\hbox{ for all $X\in\mathcal{O}_{A}$}.\]
We then in particular have that
$$m_{\s}:=\s|_{\mathcal{C}\left(\Sigma_{A}\right)} \in \Fix
\left(\mathcal{L}_{-\beta H}^{*}\right), $$
and hence  that $m_{\s}$ is a weak $(-\beta H)$--Gibbs measure.
\end{fact}
\begin{proof} Since  $\alpha_{H}^{-\i \beta}\left(S_{j}^{*}\right)=S_{j}^{*}\e^{-\beta H}$
  for all  $j\in I$, a combination of  ({\ref{CKrelation}}) and Definition 
  \ref{gauge-KMS} (2)  gives that
\begin{eqnarray*}
\s\left(X\right) &=& \sum_{j\in I}\s\left(XS_{j}S_{j}^{*}\right)
    = \sum_{j\in I}\s\left(XS_{j} \alpha_H^{\i\beta}\circ\alpha_H^{- \i\beta}(S_{j}^{*})\right)\\
&=&  \sum_{j\in I}\s\left(S_{j}^{*}\e^{-\beta H}XS_{j}\right).
\end{eqnarray*}
The assertion $\mathcal{L}_{-\beta H}^{*} m_{\s}=m_{\s} $ is an immediate consequence of (\ref{dual}).
\end{proof}
The following fact gives a generalisation of \cite[Lemma 7.5]{KerrPinzari:02}. 
In here, we have put $\{H =0\} := \{x \in \Sigma_A \with H(x) =0\}$, 
which is a closed subset of $\Sigma_A$.
\begin{fact}[\textbf{Gauge invariance}]  \label{gaugeinvariance} For an $(H,\beta)$--KMS 
state $\s$ on $\mathcal{O}_{A}$, we have
\[ \s(S_\v S_{\w}^*) = 0, \textrm{ for  all } \v,\w \in \mathcal{W}^n,\,
\v  \neq \w.\]
Furthermore, if  $\beta \in \R \setminus \{0\}$ and
 $H \geq 0$
such that $m_\sigma\left({\{H=0\}}\right) =0$, 
then
\[ \s(S_\v S_{\w}^*) = 0, \textrm{ for  all } \v, \w \in \mathcal{W},\, \v \neq \w .\]
In particular, we hence have $\s(X) =0$, for all $X$ in the closure of the vector space generated by 
$\{ S_\v S_{\w}^* \with  \v \in \mathcal{W}^m, \w\in  \mathcal{W}^n, {m \neq n} \}$.
\end{fact}
\begin{proof}  Let $\v = (\v_0 \ldots \v_{m-1})$, $\w = (\w_0  \ldots
\w_{n-1}) \in \mathcal{W}$ be given. We then have by  Fact \ref{KMS condition}
and Proposition 
\ref{rules} (6),
\begin{eqnarray*}
\s(S_\v S_{\w}^* ) &=& \s(S^*_\w \alpha_H^{\i \beta}  S_{\v})
= \s\left(  S_{\w}^*  e^{-\beta \sum_{k=0}^{m-1} H\circ \theta^k} S_{\v}  \right) \\
&=& \s\left(   e^{-\beta \sum_{k=0}^{m-1} H\circ \theta^k \circ \tau_{\w}} \1_{\theta^n[\w]} S_{\w}^*  
S_{\v}  \right).
\end{eqnarray*}
For $n=m$ Proposition \ref{rules} (2) then gives the first assertion. For $n\neq m$ note that 
Proposition \ref{rules} (2) implies that $\s(S_\v S_{\w}^* )=0$ if either $m>n$ and
$(\v_0 \ldots \v_{n-1})\neq \w$ or  $m<n$ and $(\w_0 \ldots \w_{m-1})\neq \v$. For the remaining 
statements we only consider the case $m < n$. The other cases can be dealt with in an analogous way. 
Using Facts \ref{Analytic}, \ref{betainvariance} and \ref{fact-invariance},  
we obtain for each $t \in \R$ and $\v,\w$ such that $(\w_0 \ldots \w_{m-1}) = \v$,
\begin{eqnarray*}
\s(S_\v S_{\w}^* ) &=& \s( \alpha_H^{\i t} (S_\v S_{\w}^*))
= \s\left( e^{-t \sum_{k=0}^{m-1} H\circ \theta^k} S_{\v}  S_{\w}^*
e^{t \sum_{k=0}^{n-1} H\circ \theta^k}\right) \\
&=& \s \left(  S_\v S_\w^* e^{t \sum_{k=m}^{n-1} H\circ \theta^k} \right)  
\end{eqnarray*}
Recall that by the  Cauchy--Schwartz inequality  $|\sigma(X^*Y)|^2 \leq \sigma(X^*X)\sigma(Y^*Y)$
 for all $X,Y \in \OA$  
(see \cite[Lemma 2.3.10]{BratteliRobinson:87}). Hence, by the monotone convergence theorem, we have
\begin{eqnarray*}
\left|\s(S_\v S_{\w}^* )\right|^2 & = & \left|\s 
\left(  S_\v S_\w^* e^{t \sum_{k=m}^{n-1} H\circ \theta^k} \right) \right|^2\\
&\leq&\s \left( S_\v S_{\w}^*  S_\w S_{\v}^* \right)  \sigma 
\left(e^{2t \sum_{k=m}^{n-1} H\circ \theta^k}\right)\\
&=&\s \left( S_\v S_{\w}^*  S_\w S_{\v}^* \right)  m_\sigma 
\left(e^{2t \sum_{k=m}^{n-1} H\circ \theta^k}\right) \\
&\to&\s \left( S_\v S_{\w}^*  S_\w S_{\v}^* \right)  
m_\sigma \left(\1_{\cap_{k=m}^{n-1}\{H\circ\theta^k 
= 0\}}\right), \hbox{  for }  t \to -\infty.  
\end{eqnarray*}
Since $m_\sigma \left(\1_{\cap_{k=m}^{n-1}\{H\circ\theta^k = 
0\}}\right)\leq m_\sigma \left(\1_{\{H\circ\theta^m = 0\}}\right)$,
it is now sufficient to show that $m_\sigma \left(\1_{\{H\circ\theta^m = 
0\}}\right)=0$, for all $m \in \N$.  Indeed, using   
 $\mathcal{L}_{-\beta H}^{*} m_{\s}=m_{\s}$ (see Fact 
 \ref{L*invariance}), Proposition \ref{rules} (4) and 
the assumption $m_\sigma\left(\1_{\{H=0\}}\right) =0$, 
it follows
   \begin{eqnarray*}                                                                                 
    m_\sigma \left(\1_{\{H\circ\theta^m = 0\}}\right)& = &
    m_\sigma \left(\1_{\{H = 0\}}\circ\theta^m\right)   
    = m_{\s} \left(\mathcal{L}_{-\beta H}^{m}(\1_{\{H=0\}} \circ 
    \theta^{m})\right) \\
    &=& m_{\s}\left( \sum_{v \in \mathcal{W}^{m}} S_{v}^{*}e^{-\beta 
    H} \1_{\{H=0\}}\circ \theta^{m} S_{v} \right) \\
    &=& m_{\s}\left( \sum_{v \in \mathcal{W}^{m}} S_{v}^{*}e^{-\beta 
    H} S_{v} \1_{\{H=0\}}\circ \theta^{m} \circ \tau_{v} \right) \\
    &=&  m_\sigma \left(\1_{\{H = 0\}}\mathcal{L}^m_{-\beta H}\1\right)=0.                                   
  \end{eqnarray*}   
\end{proof}

\begin{fact}[\textbf{KMS vs. Gibbs}] \label{kms_gibbs} Let $A$ be
irreducible and $\beta \in \R \setminus\{0\}$.
Then the assignment  $\Theta(\s) := m_\s:=\s|_{\mathcal{C}
\left(\Sigma_{A}\right)}$ gives rise to a well--defined surjective map
\begin{eqnarray*}
\Theta:  \{ \sigma \with  \sigma \hbox{ is a } (H, \beta)\hbox{--KMS state}  \}
&\to  & \Fix \left(\mathcal{L}_{-\beta H}^{*}\right).
\end{eqnarray*}
Furthermore, if
$H\geq0$ then the restriction of $\Theta$  to $\mathcal{S}_{0}:=\{ \sigma \with  \sigma$  
is a  $(H, \beta)$--KMS state such that 
$m_\sigma\left({\{H=0\}}\right) =0 \}$ is a bijection onto its image.
\end{fact}

\begin{proof} By the Riesz representation theorem, images under $\Theta$ are in fact 
    Borel probability measures,
    which combined with Fact \ref{conformality} gives that $\Theta$ is well--defined. 
    In order to show surjectivity
    of $\Theta$, we have to construct a state on $\mathcal{O}_A$ 
    from a given Borel measure  $\mu \in  \Fix  \left(\mathcal{L}_{-\beta H}^{*}\right)$.
For this  we define for each $v,w \in \mathcal{W}$,
\[ 
P(S_v S_w^*) := \left\{\begin{array}{l @{\quad : \quad}l} 0 & w \neq v \\ \1_{[v]} & w=v.
\end{array} \right.
\]    
This can be extended to a linear map from $\mathcal{F}_A \to \mathcal{C}(\Sigma_A)$ with $\|P(A)\|\leq \|A\|$ for each $A\in \mathcal{F}_A$ and $\|P(A)\| = \|A\|$ for $A \in \mathcal{F}_A \cap \mathcal{C}(\Sigma_A)$. Hence, $P$ admits a further extension to  a continuous and linear projection $P:\mathcal{O}_A \to \mathcal{C}(\Sigma_A)$. This allows to define the functional $\sigma_\mu$ by
\[\sigma_\mu(A):= \int P(A)d\mu \quad\hbox{ for }A\in  \mathcal{O}_A. \]
In particular, $\s_\mu $ is a positive linear
functional with $\| \s_\mu \|=1$, and hence is a state. 
In order to see that $\s_\mu$ satisfies the  KMS condition
if  $\mu \in \Fix \left(\mathcal{L}_{-\beta H}^{*}\right)$,  
it is sufficient to show that for
each $w,v,w',v'  \in \mathcal{W}$,
\begin{equation}\s_\mu(S_v S_w^* S_{v'} S_{w'}^*) = \s_\mu( S_{v'} S_{w'}^*\alpha_H^{i\beta}
    (S_v S_w^*)).\label{kms_eq} \end{equation}
For this, let $v = (v_0 \ldots v_{m-1})$, $w = (w_0 \ldots w_{n-1})$, $v' = (v'_0 \ldots v'_{p-1})$ and
$w' = (w'_0 \ldots w'_{q-1})$, for $m,n,p,q \in \N$.
Note that  by Proposition \ref{rules} (2), $ S_w^* S_{v'}= 0$ if  and only if either 
$ [w] \subset [v']$ or $ [w] \supset [v']$.
We only consider the first case,  that is $p \leq n$ and $ v'=( w_0 \ldots w_{p-1})$, 
and remark that the second case follows by exactly the same means. By Proposition 
\ref{rules} (2), we then have
\begin{eqnarray*}
S_v S_w^* S_{v'} S_{w'}^*  =   S_v S_{w_p \ldots w_{n-1}}^*  S_{w'}^*  .
\end{eqnarray*}
For $ S_v S_w^* S_{v'} S_{w'}^* \in \COA$, it immediately follows that $v = w' w_p \ldots w_{n-1}$  and that
$S_v S_w^* S_{v'} S_{w'}^* = \1_{[v]}$. Hence, $\s_\mu(S_v S_w^* S_{v'} S_{w'}^*) = \mu([v])$.
For the right hand side of  (\ref{kms_eq}) we have
\begin{eqnarray*}
    S_{v'} S_{w'}^*  \alpha^{i\beta}_H   (S_v S_w^*)
& = &  S_{v'} S_{w'}^*  e^{ -\beta H \sum_{k=0}^{m-1} H \circ \theta^k}  
S_v S_w^*  e^{\beta H \sum_{k=0}^{n-1} H \circ \theta^k}\\
&=& \1_{[w]} e^{ -\beta  \sum_{k=0}^{m-1}
 H \circ \theta^k\circ \tau_v \circ \theta^n +  \beta  \sum_{k=0}^{n-1} H \circ \theta^k}.
\end{eqnarray*}
In particular, it hence follows $ S_{v'} S_{w'}^*  \alpha^{i\beta}_H (S_v S_w^*) \in \COA$. The
latter calculation also shows  $\s_\mu(S_v S_w^* S_{v'} S_{w'}^*) \neq 0$ if and only if
$\s_\mu( S_{v'} S_{w'}^*\alpha_H^{i\beta} (S_v S_w^*)) \neq 0$. Since $ [v] =\tau_v \circ
\theta^n [w]$ and since  $d\mu\circ \theta/ d\mu = e^{\beta H}$, it follows
\begin{eqnarray*}
\s_\mu(S_v S_w^* S_{v'} S_{w'}^*)  & = & \int \1_{[v]} d \mu = 
\int \1_{[w]} \frac{d\mu \circ \tau_v \circ \theta^n}{d\mu} d \mu\\
&=& \int  \1_{[w]} \frac{d\mu \circ \tau_v}{d\mu}\circ \theta^n 
\cdot  \frac{d\mu\circ \theta^n}{d\mu} d \mu\\
&=&  \int  \1_{[w]} e^{-\beta H \sum_{k=0}^{m-1} H \circ \theta^k\circ \tau_v \circ \theta^n} 
\cdot e^{\beta H \sum_{k=0}^{n-1} H
\circ \theta^k} d \mu\\
& = & \s_\mu( S_{v'} S_{w'}^*\alpha_H^{i\beta} (S_v S_w^*)).
\end{eqnarray*}
This gives the identity in  (\ref{kms_eq}). 
To finish the proof, let
$H\geq 0$ and let $\sigma_{1}, \sigma_{2} \in \mathcal{S}_{0}$   such that 
$\Theta(\sigma_{1})=\Theta(\sigma_{2})$. Clearly, by definition we then have that 
$\s_{1}|_{\mathcal{C}\left(\Sigma_{A}\right)}= m_{\sigma_{1}}=m_{\sigma_{2}} =
\s_{2}|_{\mathcal{C}\left(\Sigma_{A}\right)}$. Furthermore,  
the gauge
invariance in Fact \ref{gaugeinvariance} gives that $\s_{1}(X) = 
\s_{2}(X)=0$, for all $X$ in the closure of the vector space generated by 
$\{ S_\v S_{\w}^* \with  \v \in \mathcal{W}^m, \w\in  \mathcal{W}^n, {m \neq n} \}$.
This shows that  the restriction of $\Theta$  to 
$\mathcal{S}_{0}$ is a bijection.
\end{proof}
\section{Radon--Nikodym representations of Cuntz--Krieger algebras}
In this section  we consider representations of Cuntz--Krieger
algebras which are
induced by  Markov fibred
systems. These representations are given by operators which act on some $L^{2}$ space
and which are determined by the Radon--Nikodym derivative of the measure 
associated with the given Markov fibred system. 
We begin with the definition of a Markov fibred system (see \cite{ADU}). Note that a 
Markov fibred system is often also referred to as  a Markov map (see e.g. \cite{Aaronson}).

\begin{defn} For a Polish space $\Omega$ and a finite Borel measure $m$ on $\Omega$, let   
$\theta:\Omega\to\Omega$ be locally invertible non--singular transformation. 
A {\em Markov partition} is a countable partition $\alpha$ of $\Omega$ for which
the following holds. $\theta:a \to \theta(a)$ is invertible, $\theta(a)$ is 
contained in the $\sigma$--algebra generated by $\alpha$ for all $a \in \alpha$, 
and the $\sigma$--algebra generated by $\{  \alpha^n :=  \bigvee_{k=0}^{n-1} 
\theta^{-k}(\alpha) \with n \in \N\}$ coincides with the Borel $\sigma$--algebra
 associated with $\Omega$, up to sets of measure zero. If $\alpha$ is a
 Markov partition, then the system  $(\Omega, m,\theta,\alpha)$ is called a\/  
 \emph{Markov fibred system}.
\end{defn}

Note that the restriction of $\theta^n$ to some arbitrary $a \in 
\alpha^n$ is clearly also always invertible and non--singular, for each $n \in \N$.
In analogy to the previous section, let $\tau_{a}$ refer to the  inverse branch of
 $\theta^n$ restricted to $a \in \alpha^n$.
Throughout, we will always assume that the Radon--Nikodym derivative
$\frac{dm\circ \theta}{dm}$ has a continuous version on the support of $m$. 
Also, if $\alpha = \{a_i \with i \in I\}$ is a finite partition then  
let $A=\left(a_{ij}\right)_{i,j\in I}$ be the finite incidence matrix
arising from $\theta$ (see e.g.  \cite{DenkerGrillenbergerSigmund:76}).
To each generator $S_{i}$ of the Cuntz--Krieger algebra
$\mathcal{O}_{A}$ we  associate an operator $s_{i}$ on  $L^{2}(\Omega,m)$
as follows.
\begin{defn}
    \label{definition:RNR}
 Let  $(\Omega, m,\theta,\alpha)$ be a Markov fibred system with associated
 incidence matrix $A=\left(a_{ij}\right)_{i,j\in I}$. For each $i\in
 I$, let $s_{i}:L^{2}(\Omega,m)\rightarrow L^{2}(\Omega,m)$ be defined by
\[
s_{i} :  f\mapsto\1_{a_{i}}\cdot{\textstyle \left(\frac{dm\circ \theta}{dm}\right)^{\frac{1}{2}}}
\cdot f\circ \theta.
\]
The $C^{*}$--operator algebra $\mathcal{R}_{A}(m)$ on the Hilbert 
space $L^{2}(\Omega,m)$  generated by the set $\{s_{i} \with i \in I\}$ will be referred to as the
{\em Radon--Nikodym representation} of  $\mathcal{O}_{A}$ induced by $(\Omega, m,\theta,\alpha)$.
\end{defn}
In order to see that a  Radon--Nikodym representation of
$\mathcal{O}_{A}$ induced by  a Markov fibred system is in fact a 
representation of the 
Cuntz--Krieger algebra $\OA$, we first make the following observations.
\begin{lem}
\label{lemma:adjoint} For each ${i}\in I$, we have that the operator $s_{i}$ is well--defined, and that
its adjoint operator $s_{i}^{*}:L^{2}(\Omega,m)\rightarrow L^{2}(\Omega,m)$ is given by
 \[
s_{i}^{*}: f\mapsto\1_{\theta(a_{i})}\cdot{\textstyle
 \left(\frac{dm\circ \tau_{a_{i}}}{dm}\right)^{\frac{1}{2}}}\cdot f\circ \tau_{a_{i}}.
 \]
 Furthermore, for each $i\in I$ and $f\in L^{2}(\Omega,m)$ we have 
 \begin{eqnarray*}
s_{i}s_{i}^{*}(f)=\1_{a_{i}}\; f\;\;\textrm{and}\;\; s_{i}^{*}s_{i}(f)=\1_{\theta(a_{i})}\; f.
\end{eqnarray*}
\end{lem}
\begin{proof}
First note that since $m$ is non-singular, we have $(dm\circ
\theta/dm)(\omega)>0$,
for $m$--almost every $\omega \in \Omega$. For each $f\in L^{2}(\Omega,m)$ and ${i}\in I$,
we then have \[
\left(s_{i}(f),s_{i}(f)\right)=\int\overline{s_{i}(f)}s_{i}(f)dm=
\int\1_{\theta(a_{i})}|f|^{2}dm<\infty.\]
 This shows that $s_{i}$ is well-defined, and also that $s_{i}^{*}s_{i}(f)=\1_{\theta(a_{i})}\cdot f$.
Next note that for $f,g\in L^{2}(\Omega,m)$,
\begin{eqnarray*}
\int\overline{f}\; s_{i}(g)\; dm & = &
\int\overline{f(\omega)}\;\1_{a_{i}}(\omega)\;{\textstyle
\left(\frac{dm\circ \theta}{dm}(\omega)
\right)^{\frac{1}{2}}}\; g\circ \theta(\omega)\; dm(\omega)\\
 & = & \int\overline{f\circ \tau_{a_{i}}(\omega)}\;\1_{a_{i}}\circ \tau_{a_{i}}(\omega)\;
 {\textstyle \left(\frac{dm\circ \theta}{dm}(\tau_{a_{i}}(\omega))\right)^{-\frac{1}{2}}}\; 
 g(\omega)\; dm(\omega)\\
 & = & \int\overline{\left(\1_{\theta(a_{i})}(\omega)\;
 {\textstyle \left(\frac{d(m\circ \tau_{a_{i}})}{dm}(\omega)\right)^{\frac{1}{2}}}\;
 f\circ \tau_{a_{i}}(\omega)\right)}\;\; g(\omega)\; dm(\omega).
 \end{eqnarray*}
This shows that $s_{i}^{*}(f)=\1_{\theta(a_{i})}\cdot{\textstyle 
\left(\frac{dm\circ \tau_{a_{i}}}{dm}\right)^{\frac{1}{2}}}
 \cdot f\circ \tau_{a_{i}}.$
The remaining part of the lemma can now be obtained as follows.
\begin{eqnarray*}
\left(s_{i}s_{i}^{*}(f)\right)(\omega) & = & \left(s_{i}\left(\1_{\theta(a_{i})}
\cdot{\textstyle \left(\frac{dm\circ \tau_{a_{i}}}{dm}\right)^{\frac{1}{2}}}
\cdot{f\circ \tau_{a_{i}}}\right)\right)(\omega)\\
 & = & \1_{a_{i}}(\omega)\cdot{\textstyle 
 \left(\frac{dm\circ \theta}{dm}(\omega)\right)^{\frac{1}{2}}}\cdot \\
 && \1_{\theta(a_{i})}(\theta(\omega))\cdot{\textstyle 
 \left(\frac{dm\circ \tau_{a_{i}}}{dm}(\theta(\omega))\right)^{\frac{1}{2}}}
 \cdot{f\circ \tau_{a_{i}}(\theta(\omega))}\\
 & = & \1_{a_{i}}(\omega)\cdot f(\omega).\end{eqnarray*}
\end{proof}
\begin{cor}
    \label{corollary:relations}
For  the generators $\{s_{i} \with i \in I\}$  of  the Radon--Nikodym
representation $\mathcal{R}_{A}(m)$ of  a
Cuntz--Krieger algebra
$\mathcal{O}_{A}$ induced by
a Markov fibred system, we have for each $ i\in I$,
     \[
    \sum_{j\in I}s_{j}s_{j}^{*}=1,\,\,\,\,\textrm{and }s_{i}^{*}s_{i}=\sum_{j\in I}a_{ij}
    s_{j}s_{j}^{*}
    .   \]
    \end{cor}
    \begin{proof}
The stated relations are immediate consequences of Definition
\ref{definition:RNR} and
Lemma \ref{lemma:adjoint}.
    \end{proof}
The following theorem summarises the results of this section.

\begin{thm} \label{RN representation}
    Let
     $\mathcal{R}_{A}(m)$ be the Radon--Nikodym representation of  the
     Cuntz--Krieger algebra
    $\mathcal{O}_{A}$ induced by a Markov fibred system
 $(\Omega,m,\theta,\alpha)$.
If the incidence matrix $A$ is irreducible, then $\mathcal{O}_{A}$
is ${*}$--isomorphic  to $\mathcal{R}_{A}(m)$.
\end{thm}

\begin{proof} Let $\rho$ refer to the
    canonical $*$--homomorphism from $\OA$ to  $\mathcal{R}_{A}(m)$ such
    that
    $\rho(S_{i})=s_{i}$, for all $i\in I$.
By Corollary  \ref{corollary:relations} we have that the generators $\left\{
s_{i}\right\} $ of $\mathcal{R}_{A}(m)$
satisfy the same type of relations as the generators $\{S_{i}\}$ of $\mathcal{O}_{A}$
(see (\ref{CKrelation})). Since $*$--homomorphisms are continuous, one immediately obtains that $\rho$ is surjective.
In order to show injectivity, note that the kernel of $\rho$ is a two-sided closed ideal in $\mathcal{O}_{A}$. Since $\mathcal{O}_{A}$ is simple, the assertion follows.  
\end{proof}

\section{Lyapunov spectra for KMS states on Cuntz--Krieger algebras associated with Kleinian groups}
In this section we apply the results of the previous sections to a particular 
class of potential functions which were essential for the 
multifractal analysis of limit sets of Kleinian groups in 
\cite{KesseboehmerStratmann:04a}. We show that this type of 
multifractal analysis
gives rise  to interesting results
concerning the existence of KMS states on Cuntz--Krieger algebras.

Recall that a Kleinian group is a discrete subgroup of the group of orientation
preserving isometries of hyperbolic $(n+1)$--space ${\mathbb{D}}^{n+1}$ (see e.g. \cite{Beardon:Intro}).
A  non-elementary Kleinian groups $G$ is called essentially free  if $G$ has a
{Poincar\'e polyhedron} $F$ (see \cite{Maskit}) with finitely many faces 
$\{ f_1, \ldots, f_{2 \mathfrak{g}}\} $ such
that if two faces $f_i$ and $f_j$ intersect inside ${\mathbb{D}}^{n+1}$, then the
 two associated generators $g_i$ and $ g_j$ of $G$ commute.
As a consequence of  Poincar\'e's theorem (see \cite{EP}), we therefore have that for an essentially free
Kleinian group there are no relations other than those arising from cusps of rank greater than 1.
Furthermore, recall that to  each essentially free Kleinian group
$G$ we can associate the following expansive coding map $T$, 
 an analogue of the Bowen--Series map in higher
dimensions (\cite{BowenSeries}, \cite{StadlbauerStratmann:04}).
For $e_i$ denoting the image of the projection of $f_i$ from some 
fixed chosen 
point in $F$ to the boundary ${\mathbb{S}}^{n}$ 
of hyperbolic space,
let $\alpha$ be the partition of $L_r(G)$ generated by $\{e_1 \cap 
L_r(G), \ldots ,e_{2\mathfrak{g}} \cap L_r(G)\}$. To each $a = e_{i_1} 
\cap \cdots \cap e_{i_k} \cap L_r(G)\in \alpha$ we then associate some arbitrary 
fixed $j(a) \in \{ i_1 , \ldots ,i_k\}$.
With $L_r(G)$ referring to the radial limit set, that is the intersection 
of $L(G)$ with the complement of the set
of parabolic fixed points of $G$, we define
\begin{eqnarray*}
T : L_r(G) \to  L_r(G),\quad
T \arrowvert_{a} := g_{j(a)} \hbox{ for } a \in \alpha.
 \end{eqnarray*}
 For further details on the construction of this map, we refer to
 \cite{StadlbauerStratmann:04}. As  shown in \cite{StadlbauerStratmann:04},
 the system $(L_r(G), \nu, T, \alpha )$ is a Markov fibred system, 
 for which $T$ is  topologically mixing, as well as conservative 
 and ergodic with respect to the canonical $T$-invariant measure $\nu$ in 
 the measure class of the Patterson measure associated with 
 $G$ (for the construction of the Patterson measure we refer to 
 \cite{Nicholls}, \cite{Patterson}, \cite{Sullivan}). It hence follows that 
 the incidence matrix $A$ arising from the symbolic
 dynamics of $T$ is irreducible. 
  Also, we clearly have that there exists a canonical coding map 
  $\pi:(\Sigma_A,\theta) \to (L(G),T)$  for which $\pi \circ \theta = T\circ \pi$, and which is 
  one--one on $L_r(G)$.
We now introduce the relevant potential functions.
Namely, let $J $ be the continuous potential function which is given by 
\[ J : \Sigma_{A} \to \R, w \mapsto \log |T'(\pi(w))|. \]
Moreover, we define for $s \in \R$,
\begin{equation}\label{I} J_{s}:=sJ+P\left(-sJ\right) \hbox{ and } I_s := J_s + \chi - \chi \circ 
\theta, \end{equation}
where $P\left(-sJ\right)$ refers to the pressure function  of $-sJ \in 
\mathcal{C}\left(\Sigma_{A}\right)$, and 
$\chi \in \mathcal{C}(\Sigma_{A})$ is determined by   
$\mathcal{L}_{-I_s} \1 =\1$.  Note that if there are no parabolic 
elements then $e^{\chi}$ is the unique eigenfunction of 
$\mathcal{L}_{-J_s}$ associated with the eigenvalue $1$. In  case there 
are parabolic elements the significance 
of $\chi$ is slightly more involved and will be given in the proof of the 
following theorem. Finally, recall that the Hausdorff dimension $\dim_{H}(\mu)$
of a Borel measure $\mu$
on $\R^{n}$ is given by
\[  \dim_{H}(\mu):= \inf \{ \dim_{H}(E): E \hbox{ is a Borel set 
with } \mu(E)>0 \}.\]
\begin{thm}\label{main} Let $\OA$ be the Cuntz--Krieger algebra associated with an 
essentially free Kleinian group $G$. Then there
is a maximal interval $(\alpha_-, \alpha_+) \subset \R_+$ and a real analytic function 
$ s: (\alpha_-, \alpha_+) \to \R$ such that for each
  $\alpha \in (\alpha_-, \alpha_+)$ there exists a 
unique   $(I_{s(\alpha)},1)$--KMS state $\s_{s(\alpha)}$ on $\OA$, 
for which we in particular have
\begin{eqnarray}\label{schottky1}
\s_{s(\alpha)}(J)  = \alpha.
\end{eqnarray}
Furthermore, $\s_{s(\alpha)}\arrowvert_{\COA}$ is a 
$\theta$--invariant weak $(-I_{s(\alpha)})$--Gibbs 
measure, and for each $\alpha \in (\alpha_-, \alpha_+)$ we have 
\begin{eqnarray}\label{schottky2}\dim_H \left(\s_{s(\alpha)}\arrowvert_{\COA} 
\, \circ \pi^{-1}\right)  = \frac{\s_{s(\alpha)}(J_{s(\alpha)} )} {\s_{s(\alpha)}(J)}.
\end{eqnarray}
Additionally, if $G$ has no parabolic elements then in the above statements the open interval
 $(\alpha_-, \alpha_+)$ can be replaced by the closed interval  
$[\alpha_-, \alpha_+]$. (For a further discussion of the boundary 
points $\alpha_{-}$ and $\alpha_{+}$, we refer to Remark (1) below).   
\end{thm}
\begin{proof}
First note that since in (\ref{schottky1}) and (\ref{schottky2})
 only
restrictions of KMS-states to $\mathcal{C}(\Sigma_{A})$ are 
considered, it is sufficient to verify (\ref{schottky1}) and
(\ref{schottky2}) for certain fixed points of the dual 
Perron--Frobenius operator associated with some suitable potential 
function. 
Next, recall that in \cite{KesseboehmerStratmann:04a} we studied 
fractal measures and Hausdorff dimensions of the
$\alpha$-level sets
\[M_{\alpha}:=\left\{ \xi \in L(G) \with \lim_{n\to\infty} \frac{1}{n} \sum_{k=0}^{n-1}
\log|T'( T^k(\xi) )|= \alpha \right\}  , \hbox{ for } \alpha \in \R .\]
Let us first discuss the case in which $G$ has no parabolic 
elements. In this 
case the map $T$ is expanding and $\log |T'|>0$, and we begin with computing the 
Lyapunov spectrum associated with the  
`homological growth rate', that is $\alpha \mapsto 
\dim_{H}(M_{\alpha})$. Note that the following  differs from the approach in 
\cite{KesseboehmerStratmann:04a}, and hence gives an alternative 
proof of the results of  \cite{KesseboehmerStratmann:04a} for the 
case in which $G$ has no parabolic elements.  Also, in the following we assume  that the reader is 
familiar with the basic results in multifractal analysis (see e.g. 
\cite{Denker}, \cite{Falconer}, \cite{Pesin}).

From a purely algebraic point of view the presumably most obvious way to
establish a KMS state on the noncommutative algebra $\OA$  
associated with $G$ is provided by the 
 measure of maximal entropy  $m$ arising from $\theta$. Clearly, the 
 topological entropy of  $(\Sigma_{A}, \theta)$ is $h_{\mathrm{top}}:=  \log(2\mathfrak{g}-1)$,
 and hence the potential giving rise to $m$ is
the constant function $\phi:\equiv\log(2\mathfrak{g}-1)$. From this we immediately 
obtain that  here the relevant  gauge action is the  
one-parameter group of $\ast$--automorphisms
$\left(\alpha_{\phi}^{t}\right)$, 
given by 
\[ \alpha_{\phi}^{t} (S_j):= e^{it \phi}  \, S_j = (2\mathfrak{g}-1)^{it} 
\, S_j \hbox{ for all } j \in \{1, \ldots, 2 \mathfrak{g}\}, t \in {\mathbb{R}}.\]
In order to determine the associated KMS state, we  
compute the temperature at which the system is at equilibrium. Observe that
for the topological pressure $P$ of the 
system we have for $s \in {\mathbb{R}}$, 
\begin{eqnarray*} P(-s \phi)&= & 
    \lim_{n \to \infty} \frac{1}{n} \log \left(2\mathfrak{g} 
    (2\mathfrak{g}-1)^{n-1} e^{-sn h_{\mathrm{top}}} 
    \right)  \\ 
&=& 
(1-s) h_{\mathrm{top}}.
\end{eqnarray*}  
This shows that  $P(-s \phi) =0$ if and only if $s=1$, and hence
the system is at equilibrium exactly for the inverse temperature 
$\beta=1$. Therefore, the KMS state  canonically 
associated with  $\phi$ is the $(\phi, 1)$--KMS 
state $\sigma_{\phi}$.  Note that by Fact \ref{kms_gibbs} we in particular 
have that 
$\sigma_{\phi}\arrowvert_{\mathcal{C}(\Sigma_{A})} =m$. Also,  one immediately
verifies that $m$ is a $\theta$-invariant $(-\phi)$--Gibbs measure for the  
H\"older continuous potential $\phi$.  Therefore,  we can apply  
standard multifractal analy\-sis, which gives that  for the level--sets 
\[ F_{\beta}:= \left\{ w \in \Sigma_{A}: \lim_{r\to 0} \frac{\log 
m\left(\pi^{-1}(B(\pi(w),r))\right)}{\log r} = 
\beta\right\}\]
we have 
\[ \dim_{H}(\pi(F_{\beta(q)})) = s(q) + q \beta(q).\]
In here, the function $s:\R \to \R$ is determined by the equation 
$P(-s(q)J -q \phi)=0$, and the function $\beta$ is given by 
$\beta(q):= - s'(q)$. Also, we have that the function given by $\beta(q) \mapsto 
\dim_{H}(\pi(F_{\beta(q)}))$ is real 
analytic 
on the image of $\beta$, which is a  closed interval 
$[\beta_{-},\beta_{+}]$.
For the $\theta$--invariant, ergodic $I_{s}$--Gibbs measure 
$m_{s(q)}$ (see (\ref{I}))   one then immediately verifies
\[ \lim_{n \to \infty} 
\frac{\sum_{k=0}^{n-1}\phi\circ\theta^k(w)}{\sum_{k=0}^{n-1} J\circ\theta^k (w)} =
\beta(q) , \hbox{  for $m_{s(q)}$--almost every } w \in \Sigma_{A}.\]
We now make the following observation, where we have set
$ \alpha(q):= \frac{h_{\mathrm{top}}}{\beta(q)}$.
\begin{eqnarray*}
    \dim_{H}(\pi(F_{\beta(q)}))  & = & \dim_{H}\left(\pi\left( \left\{ w \in \Sigma_{A}:
    \lim_{n \to \infty} {\textstyle \frac{\sum_{k=0}^{n-1}\phi\circ\theta^k(w)}{\sum_{k=0}^ {n-1}
     J\circ\theta^k (w)}} = \beta(q)  \right\} \right)\right) \\
   & = & \dim_{H}\left(\pi\left( \left\{ w \in \Sigma_{A}: \lim_{n \to  \infty}
{\textstyle \frac{ \sum_{k=0}^{n-1} J\circ\theta^k (w)}{ n}} = \alpha  (q) \right\}\right)\right).
\end{eqnarray*}

Note that this shows in particular that in the absence of parabolic 
elements the multifractal spectrum of 
the measure of maximal entropy coincides with the Hausdorff dimension
spectrum
of the homological growth rates
considered in \cite{KesseboehmerStratmann:04a}. Summarising the above 
 in terms of $\alpha$, we now have
\[  \dim_{H}(m_{s(q)} \circ \pi^{-1})  = s(q) + q \frac{h_{\mathrm{top}}}{\alpha(q)},\]
where $ s(q)$ is given by $P(-s(q) J- q \phi) 
=0$, or what is equivalent 
$P(-s(q) J)= q h_{\mathrm{top}}$. Using $P'(-s(q) J) s'(q) = 
h_{\mathrm{top}}$, we can now proceed as 
follows. First note that, by the above,
\begin{equation}\label{eq:1}  \alpha(q)= \frac{h_{\mathrm{top}}}{\beta(q)} = 
\frac{h_{\mathrm{top}}}{-s'(q)} = -P'(-s(q) J) = \int J \, d 
m_{s(q)},
\end{equation}
and that the function given by $\alpha(q) \mapsto 
\dim_{H}(m_{s(q)} \circ \pi^{-1})$ is real 
analytic 
on the image of $\alpha$, which is a closed interval 
$[\alpha_{-},\alpha_{+}]$. Furthermore, we have 
\begin{equation}\label{eq:2}
    \aligned  \dim_{H}(m_{s(q)} \circ \pi^{-1})  &=  \frac{s(q) \alpha(q) + P(-s(q) J)}{\alpha(q)}  \\
    & = \frac{-s(q) P'(-s(q) J) + P(-s(q) J)}{\alpha(q)} \\
    & =  \frac{\int \left( s(q) J + P(-s(q) J) \right) \, d m_{s(q)}}{\int  
    J \, d m_{s(q)}} \\
    & =    \frac{m_{s(q)}\left(J_{s(q)} \right)}{m_{s(q)} 
    \left( J \right) } . 
    \endaligned
\end{equation}

We now use Fact \ref{kms_gibbs} which gives that there exists an  
$\left(I_{s(\alpha)},1\right)$--KMS state
$\s_{s(\alpha)}$ on $\OA$  such that  
$\Theta(\s_{s(\alpha)}) =m_{s(\alpha)}$.
Hence, by combining this with
(\ref{eq:1}) and (\ref{eq:2}), the statements in (\ref{schottky1})  and 
(\ref{schottky2}) follow. The fact that $\s_{s(\alpha)}$ is unique can be seen as
follows.  Fact \ref{kms_gibbs} states that $\Theta$ is injective if
the underlying potential is strictly positive.  Therefore, it is
sufficient to show that $I_{s(\alpha)} >0$. This follows, since
$\mathcal{L}_{-I_{s(\alpha)} } \1 =\1$ and hence,
\begin{equation}\label{norm} 
\sum_{u \in \theta^{-1}(\{w\})}
 e^{-I_{s(\alpha)} (u)} =1, \hbox{ for all } w \in\Sigma_{A}.
\end{equation}
Finally, note that since  ${-I_{s(\alpha)} }$ is H\"older
continuous, the Perron-Frobenius-Ruelle Theorem (see e.g. \cite{Bowen}) implies
that $\Fix \left(\mathcal{L}^*_{-I_{s(\alpha)} }
\right)$  is a singleton. This finishes the proof if $G$ has no 
parabolic elements.

We now consider the parabolic situation. Hence,  let  $G$ be an essentially
 free Kleinian group with parabolic elements. 
It is well--known that in this case the limit set $L(G)$ can be 
written as the disjoint union of $L_r(G)$ and the countable set of fixed points 
of the parabolic transformations in $G$ (see \cite{BM}). In contrast to 
the previous case, $T$ is now expansive and the function $\log |T'|$
is equal to zero precisely on the fixed points of the parabolic generators of $G$. 
In addition to the coding by $\Sigma_A$, there is an alternative coding which
 is provided by the method of inducing. That is, for 
a subset $B$ of $\Sigma_{A}$ we obtain the 
induced map 
$\widetilde{\theta}: B \to B, w \mapsto \theta^{N(w)} (w)$, where $ 
N: \Sigma_{A} \to 
\N \cup \{\infty\}, w \mapsto \inf\{n \in \N \with \theta^{n}(w)\in B \}$ 
(for further details we refer to the Appendix). We always 
assume that $\pi(B)$ is bounded away from the set of fixed points of the 
parabolic generators of $G$. This guarantees that $\widetilde{\theta}$
is an expanding Markov fibred system 
$(B,\widetilde{\nu},\widetilde{\theta}, \widetilde{\alpha})$ with
 respect to a countable partition of 
$B$, where $\widetilde{\nu}$ refers to the pull--back under $\pi^{-1}$ 
of the invariant version of the restriction of the Patterson measure to 
$\pi(B)$.  For further details 
on the construction of this system we refer to \cite{StadlbauerStratmann:04}.  
The canonical potential function $\widetilde{J}$ for the induced system 
is given for $w \in B$ by
\[ \widetilde{J}(w) :=\sum_{k=0}^{N(w)-1} J(\theta^k(w)). \]
As shown in  \cite{KesseboehmerStratmann:04a}, in this parabolic situation there is a maximal interval 
$(\alpha_-, \alpha_+)\subset \R_+$ and a real analytic function 
$ s: (\alpha_-, \alpha_+) \to \R$ such that
for each  $\alpha \in (\alpha_-, \alpha_+)$ there exists  a unique  
$(-\widetilde{J}_{s(\alpha)})$--Gibbs measure 
$\widetilde{\mu}_{s(\alpha)}$, where $\widetilde{J}_{s(\alpha)} :=
 {s(\alpha)}\widetilde{J} + P(-{s(\alpha)}J)N$. In particular, 
$\widetilde{\mu}_{s(\alpha)}( B \cap \pi^{-1}(M_\alpha))=1$ and
\begin{equation}\label{eq:3}
      \dim_{H}(\widetilde{m}_{s(\alpha)} \circ \pi^{-1})  = \dim_{H}(M_\alpha) = 
      \frac{\int   \widetilde{J}_{s(\alpha)} d \widetilde{m}_{s(\alpha)}}{\int 
\widetilde{J} d \widetilde{m}_{s(\alpha)}}. 
      \end{equation}
Moreover, with $\widetilde{h}$ referring to the eigenfunction for the eigenvalue 
1 of the Perron--Frobenius operator 
$\widetilde{\mathcal{L}}_{-\widetilde{J}_{s(\alpha)}}$ of the induced system, let
\[   
\widetilde{I}_{s(\alpha)} := {s(\alpha)}\widetilde{J} + 
P(-{s(\alpha)}J)N + \log\widetilde{h} - \log\widetilde{h} \circ \widetilde{T}. 
\]
Then there exists a unique $\widetilde{\theta}$--invariant $(-\widetilde{I}_{s(\alpha)})$--Gibbs 
measure $\widetilde{m}_{s(\alpha)}$ in the measure class of $\widetilde{\mu}_{s(\alpha)}$. 
Since $\widetilde{m}_{s(\alpha)}$ is $\widetilde{\theta}$--invariant  we 
obtain a $\theta$--invariant measure $m_{s(\alpha)}$ by Kac's formula (see Appendix). 
 Clearly,  the measures $\widetilde{\mu}_{s(\alpha)}$, 
$\widetilde{m}_{s(\alpha)}$  and $m_{s(\alpha)}$ are all contained in the same 
measure class. 
By setting $H= J_{s(\alpha)}$ in
Corollary \ref{coboundary}, it follows that there exists a continuous
function $\chi : \Sigma_{A}\to \R$ such that for the
Radon--Nikodym derivative of $m_{s(\alpha)}$ we have (see Corollary
\ref{coboundary} and the remark thereafter)
\[  \frac{d{m_{s(\alpha)}}\circ{\theta}}{d{m_{s(\alpha)}}} = e^{J_{s(\alpha)}+
\chi - \chi \circ \theta },\]
or equivalently $m_{s(\alpha)}\in 
\Fix\left(\mathcal{L}^*_{-I_{s(\alpha)}}\right)$, where
$ I_{s(\alpha)}:= s(\alpha) J+P(-s(\alpha) J)+ \chi - \chi \circ  \theta$.
Therefore, \cite{Kesseboehmer:01} implies that $m_{s(\alpha)}$ is a 
$\theta$--invariant weak $(- I_{s(\alpha)})$--Gibbs measure. 
Next we show uniqueness of $m_{s(\alpha)}$. First note that
\cite[Proof of Theorem 1.2]{KesseboehmerStratmann:04a} shows that 
equilibrium measures for $-I_{s(\alpha)}$
are mapped to equilibrium measures for $-\widetilde{I}_{s(\alpha)}$ 
by inducing. If we restrict to measures having full
measure on $\bigcup_{i=0}^\infty \theta^{-i} (B)$, then the inverse of 
this mapping is given by Kac's formula.
Since by \cite[Theorem 2.2.9]{MU} the set of equilibrium measures for 
$-\widetilde{I}_{s(\alpha)}$ is a singelton, it is 
now sufficient to show that for every ergodic 
equilibrium measure $\nu$ for the potential $-I_{s(\alpha)}$ we have 
$\nu(\bigcup _{i=0}^\infty \theta^{-i} (B))=1$. Indeed, since the complement of 
$\bigcup _{i=0}^\infty \theta^{-i} B$ corresponds to the countable set of
 parabolic fixed points of $G$, we have that
$\nu(\Sigma_A \setminus \bigcup _{i=0}^\infty \theta^{-i} B)=1$ implies 
$h_\nu=0=\int -s(\alpha)I\,d\nu=0$.
Hence, since $\nu$ is an equilibrium state, it follows $P(-I_{s(\alpha)})=0$.
This contradicts the fact that $P(-I_s)>0$ for
all $s<\delta$, and therefore gives the uniqueness of $m_{s(\alpha)}$.

In order to verify the statements in (\ref{schottky1}) and (\ref{schottky2}), observe that by 
construction of $m_{s(\alpha)} $ we have 
 \begin{equation}\label{eq:4}
\int J d m_{s(\alpha)}  =  \frac{1}{\widetilde{m}_{s(\alpha)} (N)} \int 
\widetilde{J} d \widetilde{m}_{s(\alpha)} \end{equation}
and
\begin{equation}\label{eq:5}
 \int J_{s(\alpha)} dm_{s(\alpha)}  =  \frac{1}{\widetilde{m}_{s(\alpha)} (N)} \int  
 \widetilde{J}_{s(\alpha)} d \widetilde{m}_{s(\alpha)}.
\end{equation}
Since $I_{s(\alpha)}$ is strictly positive for $\alpha \in 
(\alpha_{-},\alpha_{+})$,  we can now apply Fact \ref{kms_gibbs} as in 
the previous case. It follows  
that there  exists a  $\left(I_{s(\alpha)},1\right)$--KMS 
state $\s_{s(\alpha)}$ such that
$\s_{s(\alpha)}\arrowvert_{\COA}=m_{s(\alpha)} $.  Hence, by combining this 
with (\ref{eq:3}), (\ref{eq:4}) and (\ref{eq:5}), the assertions in (\ref{schottky1}) 
and (\ref{schottky2}) follow. 
\end{proof}

By combining Theorem \ref{RN representation} and Theorem \ref{main},
we immediately obtain the following result. For further details on
the relations between KMS states and vector states we refer to 
\cite{BratteliRobinson:87}.

\begin{cor} Let $\OA$ be the Cuntz--Krieger algebra associated
    with an essentially free Kleinian group $G$, and let
    $ s: (\alpha_-, \alpha_+) \to \R$ and $m_{s(a)}$ be
    given by the previous theorem. Then there exists a real
    analytic family $\{ \mathcal{R}_A(m_{s(a)}) \with a \in (\alpha_-, 
    \alpha_+)\}$ of  faithful Radon--Nikodym representations
    $ \mathcal{R}_A(m_{s(a)})$ induced by the Markov 
    fibred system $(L(G),m_{s(a)}\circ \pi^{-1},T,\alpha)$. In 
    particular, for each $a \in (\alpha_{-},\alpha_{+})$ we have 
    that the 
    $(I_{s(a)},1)$--KMS state $\s_{s(a)}$ in the previous 
    theorem is a vector state which is given by
    \[ \sigma_{s(a)} (X) = \left(\1 , X (\1)\right)_{s(a)}.\]
    In here, the inner product $(\cdot , \cdot)_{s(a)}$ refers to the
    inner product 
    on the Hilbert space 
    $\left(L^{2}(L(G)), m_{s(a)}\circ \pi^{-1}\right)$.
\end{cor}

{\bf Remarks.} \\ 
(1) \,  Note that  $\alpha_{+}$ is  always given by
  $\alpha_+ = \lim_{s \to - \infty} P(-sJ)/(-s)$.  Similar, if $G$ has no 
  parabolic elements then  $\alpha_- =
  \lim_{s \to  \infty} P(-sJ)/(-s)$. Whereas, if $G$ has parabolic elements
  then  $\alpha_-= \lim_{s \nearrow \delta} 
  P(-sJ)/(\delta - s)$, where $\delta=\delta(G)$ refers to the exponent 
  of convergence of $G$.  Here, the parabolic case has to be 
  treated with extra care, since as we have  shown 
  in \cite{KesseboehmerStratmann:04a} 
  in this situation a phase transition 
  can occur at $\delta$ (see also \cite{StadlbauerStratmann:04},
  \cite{KS1}, \cite{KS2}). 
  More precisely,  we have the following scenario, where $k_{\mathrm{max}}$
  refers to 
  the maximal possible rank of the parabolic elements in $G$.
 For $\delta \leq (k_{\mathrm{max}}+1)/2$ we have that $\alpha_{-}=0$. In this
 situation one
 immediately verifies that $\widetilde{m}_{s(0)}(N)=\infty$, and 
 consequently  Kac's  formula is not applicable. However, 
 we still obtain a $\theta$--invariant probability measure 
 $m_{s(0)}$ 
 as the weak limit of a sequence $\left(m_{s(\alpha_{n})}\right)$,
 for $\alpha_{n}$ tending to $0$ from above. Note that the 
 measure $m_{s(0)}$ has to be purely atomic.
 On the other hand, if $\delta > (k_{\mathrm{max}}+1)/2$, then 
 $\alpha_{-}>0$ and $\widetilde{m}_{s(0)}(N)<\infty$.
 Hence, we can argue as in  the proof of Theorem \ref{main}
 to obtain that in this case the boundary point $\alpha_{-}$
 can be included in the statement of Theorem \ref{main}.

  (2) \,
  In order to see that the assignment $q \mapsto s(q)$ gives rise to a strictly convex
  function, one can argue as follows.  We only consider the 
  non--parabolic case, and refer to \cite{KesseboehmerStratmann:04a} for the parabolic 
  situation. Using the notation in the proof of 
  Theorem \ref{main},  we have for the second derivative of $s$ 
   (see e.g. \cite[p. 237]{Denker}), 
  \[ s''(q) = \frac{D_{q}(s'(q) J - \phi)}{\int \log |T'| \, d m_{q}},\]
  where $D_{q}$ refers to the asymptotic covariance given for 
  a H\"older continuous function $f$ on $\Sigma_{A}$ by 
  \[ D_{q}(f) := \sum_{k=0}^{\infty} \left( \int f \cdot f 
  \circ \theta^{k} \, d m_{q} - \left(\int f  \,  d m_{q} \right)^{2}  \right) .\]
  Therefore, the function $s$ is strictly convex if and 
  only if $s(q) J + q \phi$ is not cohomologous to a constant.
  In order to see that the latter does in fact hold, one can argue  similar
 as  in
 the proofs of the `dynamical rigidity theorems'  of \cite{B4}, 
 \cite{MU} and \cite{S4}. Namely, the assumption that 
  $s(q) J + q \phi$ is cohomologous to a constant  is equi\-va\-lent to
  the statement that there exists a constant $R$ such that  for all 
  $n \in \N$ and $w \in \Sigma_{A}$ for which $\theta^{n}(w) = w$
  (see e.g. \cite[Theorem 2.2.7]{MU}),
  \[ \sum_{k=0}^{n-1} \left( s(q) J(\theta^{k}(w)) +q \phi(\theta^{k}(w)) 
  \right) = n R.\]
  One immediately verifies that the latter identity is equivalent to
  \[  s(q)  \log |\left(T^{n}\right)' (\pi(w)) |
  = n ( R- q h_{\mathrm{top}}) .\]
 This shows that the periodic points of period $n$ in $L(G)$ must all
 have equal multipliers, which is clearly absurd for the conformal 
 system given by $G$. (Note that in here the constant $R$ is in fact 
 given by $R= P(q \phi )- P(-s(q) J)$).

(3) \,   Finally, we remark that if $G$ has no parabolic elements then  
 the above analysis gives rise to the following estimate of
the asymptotic
growth  rate  of the word metric in $G$ for generic elements of $L(G)$.
For this note that  the Lyapunov spectrum of
 the measure of maximal entropy  attains its maximum precisely
 at the exponent of convergence $\delta=\delta(G)$. In the notation used in the 
proof of Theorem \ref{main}, we then have $s(q)=\delta$, and 
hence since $P(\delta J)=0$
and $P(-s(q) J)= q h_{\mathrm{top}}$, it follows that $q=0$.  This implies that
\[ \beta(0)=   \frac{ h_{\mathrm{top}} }{\int 
\log |T'| \, d \nu},\]
where $\nu$ refers to the invariant version of the Patterson 
 measure $\mu$ constructed with respect to the origin 
in ${\mathbb {D}}^{n+1}$. Now,  let $\xi_{t}$ 
refer to the unique point on the ray from the origin to $\xi \in 
{\mathbb{S}}^{n}$
at hyperbolic distance $t$ to the origin. Also, let $[\xi_{t}]$ denote 
the word length of $g$, for $g \in G$  determined by 
$\xi_{t} \in g(F)$. Then the arguments in the proof of Theorem 
\ref{main}
immediately imply that  for $\mu$--almost every $\xi \in 
L(G)$ we have 
\[ \lim_{t\to \infty} \frac{t}{[\xi_{t}] } = \frac{h_{\mathrm{top}} }{\beta(0)}=
\int 
\log |T'| \, d \nu.\]

\section{Appendix}

In this appendix we give a refinement of a formula of 
 Kac in the context of Markov fibred systems. We obtain 
 explicit 
 formulae  which allow to compute the Radon--Nikodym derivative of a 
 $\theta$--invariant 
measure
on the whole system $(\Omega, \theta)$  in terms of the
 Radon--Nikodym derivative  of a $\widetilde{\theta}$--invariant 
 measure on an   
 induced system $(B,\widetilde{\theta})$.

More precisely, let $(\Omega,m,\theta,\alpha)$ be a conservative and ergodic Markov fibred system with respect 
to the finite partition $\alpha$ of $\Omega$. 
Furthermore, let $B \subset \Omega$ be measurable with respect to the $\sigma$--algebra generated
by $\alpha^n :=  \bigvee_{k=0}^{n-1} \theta^{-k}(\alpha)$, for some $n \in 
\N$. Also, let $\widetilde{\theta}$ refer to the \emph{induced transformation}  given by 
\[ \widetilde{\theta} :B \to B, \omega \mapsto \theta^{N(\omega)} (\omega),\]
where
\[ N: \Omega \to \N \cup \{\infty\}, \omega\mapsto \inf\{n \in \N \with \theta^{n}(\omega)\in B \}.\]
It is well--known that the induced system $(B,\widetilde{m},\widetilde{\theta}, \widetilde{\alpha}) $ is again a 
conservative and ergodic Markov fibred system, where $\widetilde{\alpha}$ denotes the associated countable 
partition which can be finite or infinite, and $\widetilde{m}:=m\arrowvert_B$ (see e.g. \cite{Aaronson}). 
The inverse branches of $\widetilde{\theta}$ will be denoted by $\widetilde{\tau}_a$, for 
$  a \in \widetilde{\alpha}$. 

Recall that  $\theta$--invariant measures  and $\widetilde{\theta}$--invariant measures are related as follows. 
If $\nu$ is a given $\theta$--invariant measure then we obtain a  $\widetilde{\theta}$--invariant measure  by 
restricting $\nu$ to $B$. Conversely, if $\widetilde{\nu}$ is a given
 $\widetilde{\theta}$--invariant measure such that $\widetilde{\nu}(N)<\infty$ one obtains a $\theta$--invariant 
 probability measure $\nu$ by the following  formula of Kac (see 
 \cite{Kac}). Namely,  for $ \sum_{k=0}^{N-1} f
\circ \theta^{k} \in L^1(\widetilde{\nu})$ we have
\[ \int f d\nu = \frac{1}{\widetilde{\nu}(N)}  \int_B \sum_{k=0}^{N(\omega)-1} f
\circ \theta^{k}(\omega) \, d\widetilde{\nu}(\omega).
\]

We now investigate for this situation in which way the two associated Radon--Nikodym derivatives
 ${d\nu\circ \theta}/{d\nu}$ 
and ${d\widetilde{\nu}\circ \widetilde{\theta}}/{d\widetilde{\nu}}$ are related.
One direction is immediately given by the chain rule. Namely, for a given $\theta$--invariant 
measure $\nu$ we have
\[  \log \frac{d\widetilde{\nu}\circ \widetilde{\theta}}{d\widetilde{\nu}} (\omega) = 
\sum_{k=0}^{N(\omega)-1}  \log\left(\frac{d\nu\circ \theta}{d\nu}(\theta^k(\omega)) 
\right).\]
The converse direction is slightly more delicate and will be subject
of the  following proposition.
We remark that it might be that this statement is known to
experts
in this area, 
however we were
unable to find it in the literature and hence decided to include the 
proof. 
We require the following notation. Let $D_n := \{\omega \in \Omega \with N(\omega) =n \}$, and 
put $N(A):=n$ if $A \subset D_n$ for some $n \in  \N$. Also, for $ A \subset 
D_n$ such that for some $b \in \widetilde{\alpha}$ we have that  
either $A \subset  B \cap b$ or $A \subset \Omega \setminus B$ and  $\theta^n(A)\subset b$,  
 we define
\[
 \mathcal{Z}(A)  := \left\{
 \begin{array}{l @{\quad:\quad}l }
     \{  a \in   \widetilde{\alpha} \with A \subset \theta^{N(a)}(a) \}   &   
     A \subset B \\
  \{  a \in   \widetilde{\alpha} \with N(a)>N(A), A \subset \theta^{N(a)-N(A)}(a) \}      & 
   A \subset \Omega \setminus B.
\end{array}
   \right.
\]
Furthermore, we put $ \mathcal{Z}(\omega) :=   \mathcal{Z}(b)$ if 
either 
$\omega \in b \in \widetilde{\alpha}$ or 
 $\omega \in \Omega \setminus B$ such that $\omega \in b \in \alpha^{n}$, 
 where $\theta^n(b) \in \widetilde{\alpha}$ for some $n \in \N$. Note 
 that in the first case  the set $\{\widetilde{\tau}_a \with a \in  
 \mathcal{Z}(\omega)\}$ represents the set of  inverse branches of $\widetilde{\theta}$ at 
 $\omega$, whereas in the second case the set $\{\widetilde{\tau}_a \with 
 a \in  \mathcal{Z}(\omega)\}$ refers to  the set of  inverse 
 branches of $\widetilde{\theta}$ at $\theta^{N(\omega)}(\omega)$ 
 with the additional property  that  $\omega \in \{ \theta^k(\widetilde{\tau}_a(\omega)) \with 
 1\leq  k < N(\widetilde{\tau}_a(\omega)) \}$, for each $a \in  \mathcal{Z}(\omega)$. 
 Hence,  for $\omega \notin B$ we in particular have
\[   \{\widetilde{\tau}_a(\theta^{N(\omega)}(\omega)) \with a \in  \mathcal{Z}(\omega) \} = 
\bigcup_{l = N(\omega)+1}^{\infty} \theta^{-(l-N(\omega))}(\omega) \cap D_l \cap B.\]  

\begin{prop}\label{prop:lifted_derivative} Let $(\Omega,m,\theta,\alpha)$ be a conservative and ergodic Markov 
    fibred system, and let $(B,\widetilde{m},\widetilde{\theta}, \widetilde{\alpha})$ be the induced system as 
    introduced above. If $\widetilde{\nu}$ is a given $\widetilde{\theta}$--invariant measure 
    which is absolutely 
    continuous with respect to $\widetilde{m}$, then the following holds for the $\theta$--invariant measure $\nu$ 
    obtained through Kac's formula. 
\begin{enumerate}
  \item For  $\nu$--almost all $\omega \in B$, we have
   \[ 
    \frac{d\nu\circ \theta}{d\nu}(\omega) =    
   \left( \sum_{a \in \mathcal{Z}(\theta\omega)} \frac{d\widetilde{\nu} \circ
    \widetilde{\tau}_a}{d\widetilde{\nu}} (\theta^{N(\omega)} (\omega))  \right)
    \cdot \frac{d\widetilde{\nu}\circ \widetilde{\theta}}{d\widetilde{\nu}}(\omega).
    \]
  \item For  $\nu$--almost all $\omega \in \Omega \setminus B$, we have
  \[ 
  \frac{d\nu\circ \theta}{d\nu}(\omega) \!=   \! \left(  \sum_{a \in \mathcal{Z}(\theta\omega)} 
  \frac{d\widetilde{\nu} \circ \widetilde{\tau}_a}{d\widetilde{\nu}} (\theta^{N(\omega)} (\omega))  
  \right)\!\! \left/ \!\!\left( \sum_{a \in \mathcal{Z}(\omega)}
   \frac{d\widetilde{\nu} \circ \widetilde{\tau}_a}{d\widetilde{\nu}} (\theta^{N(\omega)} (\omega))  \right)  
   \right.  \!\!.
   \]  
\end{enumerate}
\end{prop}
\begin{proof}
 First note that the infinite sums in (1) and (2) do converge. This follows since the 
 $\widetilde{\theta}$--invariance of $\widetilde{\nu}$ implies, for $\nu$--almost all $\omega \in B$,
\begin{equation} \label{transfer}
\sum_{a\in \mathcal{Z}(\omega)} \frac{d\widetilde{\nu} 
    \circ \widetilde{\tau}_a}{d\widetilde{\nu}} (\omega)= \sum_{a \in \widetilde{\alpha}} 
    \frac{d\widetilde{\nu} \circ \widetilde{\tau}_a}{d\widetilde{\nu}} (\omega) =1.\end{equation}   
    Let $A \subset D_n \setminus B$ such that $\theta^n(A) \subset a$, 
    for some $n\in \N$ and $a \in \widetilde{\alpha}$.  We then have
\begin{eqnarray*}
\nu(A)& = &  \frac{1}{\widetilde{\nu} (N)} \int_B \sum_{k=0}^{N(\omega)-1} \1_A
\circ \theta^{k}(\omega) \, d\widetilde{\nu}(\omega) \\
  &=& \frac{1}{\widetilde{\nu} (N)}  \sum_{l=1}^\infty   \int_{B\cap D_l} \sum_{k=0}^{l-1} \1_A
\circ \theta^{k}(\omega) \, d\widetilde{\nu}(\omega) \\
&=& \frac{1}{\widetilde{\nu} (N)}\sum_{l=1}^\infty  \sum_{k=0}^{l-1} 
 \widetilde{\nu} ( B  \cap D_l \cap \theta^{-k}(A)) 
\\ &=& \frac{1}{\widetilde{\nu} (N)} \sum_{l=n+1}^\infty  
\widetilde{\nu} ( B  \cap D_l \cap \theta^{-(l-n)}(A))\\
&=&  \frac{1}{\widetilde{\nu} (N)}  \sum_{a \in \mathcal{Z}(A)} 
\widetilde{\nu}(\widetilde{\tau}_a (\theta^n A)).
\end{eqnarray*}
Hence,  we have for $ \omega \in A$, 
\begin{equation} \label{para}  
 \frac{d\nu}{d\nu\circ \theta^n}(\omega) = 
 \sum_{a \in \mathcal{Z}(\omega)} \frac{d\widetilde{\nu} \circ \widetilde{\tau}_a}{d\widetilde{\nu}} 
 (\theta^n (\omega))  .
 \end{equation}
Therefore,  $ \nu(\theta A) =   \left(\sum_{a \in 
\mathcal{Z}(\theta A)} \widetilde{\nu}(\widetilde{\tau}_a (\theta^n 
A))\right) / \widetilde{\nu} (N)$ and  
\begin{equation}\label{para>1}  \frac{d\nu \circ \theta}{d\nu\circ \theta^n}(\omega) =  
\sum_{a \in \mathcal{Z}(\theta\omega)} \frac{d\widetilde{\nu} \circ \widetilde{\tau}_a}{d\widetilde{\nu}}
 (\theta^n (\omega))  . \end{equation}
Combining (\ref{para}) and (\ref{para>1}), the assertion follows in the 
case in which $\omega \in  D_n \setminus B$ for some $n >1$. 
The case $\omega \in  D_1 \setminus B$ is an immediate consequence of 
(\ref{transfer}) and (\ref{para}).
This proves the assertion in (2). The proof of (1) is now an immediate consequence of  (\ref{para>1}).
 Namely, for each $ \omega \in D_n \cap B$ with $n>1$,
\begin{eqnarray*}
  \frac{d\nu\circ \theta}{d\nu}(\omega)& = & \frac{d\nu \circ \theta}{d\nu\circ \theta^n}(\omega) 
  \left/  \frac{d\nu }{d\nu\circ \theta^n}(\omega) \right. =  \frac{d\nu \circ \theta}{d\nu\circ
  \theta^n}(\omega) \cdot \frac{d\widetilde{\nu}\circ \widetilde{\theta}}{d\widetilde{\nu}}(\omega)\\
 & = &  \left( \sum_{a \in \mathcal{Z}(\theta\omega)} \frac{d\widetilde{\nu} \circ \widetilde{\tau}_a}
 {d\widetilde{\nu}} (\theta^n (\omega))  \right)
    \cdot \frac{d\widetilde{\nu}\circ \widetilde{\theta}}{d\widetilde{\nu}}(\omega).
\end{eqnarray*}
If $ \omega \in D_1 \cap B$, then we have similar to the previous case 
that the statement is an immediate consequence of 
(\ref{transfer}).
\end{proof}

For the following we define, for $\omega \in \Omega$, 
\begin{eqnarray*}
Z(\omega) &:=& \{\eta \in B \with 
\theta^{k}(\eta)=\omega  \hbox{ for some }  1 \leq k \leq N(\eta). \}
\end{eqnarray*} 

\begin{cor} \label{coboundary} In the situation of the previous proposition assume that there exist measurable
    functions $H: \Omega \to \R$  and $\widetilde{\chi}: B \to \R$ such that for almost all $\omega \in B$,
\[ 
\frac{d\widetilde{\nu}\circ \widetilde{\theta}}{d\widetilde{\nu}} (\omega) =
 e^{\left( \sum_{k=0}^{N(\omega)-1}
H \circ \theta^k(\omega) \right)+  \log \widetilde{\chi}(\omega)- \log 
\widetilde{\chi}(\widetilde{\theta} (\omega))}.
 \]
Then there exists a function $\chi : \Omega \to \R$ such that for almost all $\omega \in \Omega$,
\[  \frac{d{\nu}\circ{\theta}}{d{\nu}}(\omega) = e^{H(\omega) + \log\chi(\omega) - 
 \log\chi(\theta(\omega))  } .\]
In here, the function $\chi$ is given by
\[  \chi(\omega) :=  \left\{ 
\begin{array}{l @{\,:\,} l }
     \widetilde{\chi}(\omega)& \omega \in B   \\
       \left( \sum_{\eta \in {Z}(\omega)} e^{ -\sum_{k=0}^{N(\eta)- N(\omega) -1} H\circ \theta^k (\eta) 
      - \log \widetilde{\chi}(\eta)}\right)^{-1}     & \omega \notin B.    
\end{array}
  \right. \]
\end{cor}

\begin{proof} Note
    that $\theta^{N(\eta)-N(\omega)}(\eta) = \omega$, for each $\omega 
    \in \Omega, \eta \in  {Z}(\omega)$. We hence have for $\omega \notin B$,
\begin{eqnarray*}
  &  &  \sum_{a \in \mathcal{Z}(\omega)}  \frac{d\widetilde{\nu}\circ \widetilde{\tau}_a}{d\widetilde{\nu}}(\theta^{N(\omega)}(\omega))
  =   \sum_{\eta \in Z(\omega)}  
  e^{-\sum_{k=0}^{N(\eta)-1} H\circ \theta^k(\eta) - \log\widetilde{\chi}(\eta) + \log\widetilde{\chi}(\widetilde{\theta}(\eta))} \\
  & = &  \left( \sum_{\eta \in {Z}(\omega)} 
  e^{-\sum_{k=0}^{N(\eta) - N(\omega) -1} H\circ \theta^k(\eta) - \log \widetilde{\chi}(\eta)} \right)
  e^{-\sum_{k=0}^{N(\omega) -1} H\circ \theta^k(\omega) +\log\widetilde{\chi}(\theta^{N(\omega)}(\omega))} \\
  &=& e^{- \log\chi(\omega)} \ e^{-\sum_{k=0}^{N(\omega) -1} H\circ \theta^k(\omega) +\log\widetilde{\chi}(\theta^{N(\omega)}(\omega))}.
  \end{eqnarray*}
Combining this with (\ref{transfer}) and Proposition \ref{prop:lifted_derivative}, 
the assertion follows. \end{proof}

{\bf Remark.} 
Note that one immediately verifies, using (\ref{transfer}), that the definition 
of $\chi$ in Corollary  \ref{coboundary}  
can be rewritten so that  
only finite sums are involved. Therefore, it follows that 
 $\chi$ is continuous whenever both $H$ and $\widetilde {\chi}$  are continuous.


\end{document}